\documentclass[journal=jacsat,manuscript=article]{achemso}

\usepackage{chemformula} 
\usepackage[T1]{fontenc} 

\usepackage{epsfig,amsfonts,color}
\usepackage[final]{pdfpages}
\usepackage{blkarray}
\usepackage{bm}
\usepackage{algorithm}
\usepackage{algpseudocode}
\usepackage{multirow}
\usepackage{graphicx}
\usepackage{amsmath}
\usepackage{algorithm}
\usepackage{booktabs} 
\usepackage{longtable}
\usepackage{subcaption}

\newcommand{\mn}{\mathcal{N}}
\newcommand{\me}{\mathcal{E}}
\newcommand{\mg}{\mathcal{G}}

\newcommand{\bx}{\boldsymbol{x}}

\newcommand{\be}{\begin{equation}}
\newcommand{\ee}{\end{equation}}
\newcommand{\bea}{\begin{eqnarray}}
\newcommand{\eea}{\end{eqnarray}}

\newcommand{\bvec}{\left(\begin{array}{c}}
\newcommand{\evec}{\end{array}\right)}
\newcommand{\bsub}{\begin{subequations}}
\newcommand{\esub}{\end{subequations}}

\author{Parth Brahmbhatt}
\affiliation[University of Wisconsin-Madison]
{Department of Chemical and Biological Engineering, University of Wisconsin-Madison, Madison, WI 53706, USA}

\author{David L. Cole}
\affiliation[University of Wisconsin-Madison]
{Department of Chemical and Biological Engineering, University of Wisconsin-Madison, Madison, WI 53706, USA}

\author{Victor M. Zavala}
\affiliation[University of Wisconsin-Madison]
{Department of Chemical and Biological Engineering, University of Wisconsin-Madison, Madison, WI 53706, USA}
\alsoaffiliation[]{Mathematics and Computer Science Division, Argonne National Laboratory, Lemont, IL 60439, USA}

\author{Styliani Avraamidou}
\affiliation[University of Wisconsin-Madison]
{Department of Chemical and Biological Engineering, University of Wisconsin-Madison, Madison, WI 53706, USA}
\email{avraamidou@wisc.edu}

\title[Title]
  {Benders Decomposition using Graph Modeling and Multi-Parametric Programming}

\begin{document}


{\em We dedicate this contribution to the memory of Dr. Pedro Castro. Pedro was a brilliant researcher that made exceptional contributions to the modeling, solution, and applications of optimization.} 

\begin{abstract}
Benders decomposition is a widely used method for solving large and structured optimization problems, but its performance is affected by the repeated solution of subproblems. We propose a flexible and modular algorithmic framework for accelerating Benders decomposition. Specifically, we express the problem structure by using a graph-theoretic modeling abstraction in which nodes represent optimization subproblems and edges represent connectivity between subproblems. A key innovation of our approach is that we embed multi-parametric programming (mp) surrogates for node subproblems. The use of mp surrogates allows us to replace subproblem solves with fast look-ups and function evaluations for primal and dual variables during the iterative Benders process. We formally show the equivalence between classical Benders cuts and those derived from the mp solution. We implement our framework in the open-source {\tt PlasmoBenders.jl} software package. To demonstrate the capabilities of the proposed framework, we apply it to a two-stage stochastic programming problem, which aims to make optimal capacity expansion decisions under market uncertainty. We evaluate both single-cut and multi-cut variants of Benders decomposition and show that the use of mp surrogates achieves substantial speedups in subproblem solve time, while preserving the convergence guarantees of Benders decomposition. We highlight advantages in solution analysis and interpretability that is enabled by mp critical region tracking; specifically, we show that these reveal how decisions evolve geometrically across the Benders search. Our results aim to demonstrate that combining surrogate modeling  with graph modeling offers a promising and extensible foundation for structure-exploiting decomposition. In addition, by decomposing the problem into more tractable subproblems, the proposed approach also aims to overcome scalability issues of mp. Finally, the use of mp surrogates provides a unifying and modular optimization framework that enables the representation of heterogeneous node subproblems as modeling objects with a homogeneous structure.
\end{abstract}

\section{Introduction}
\label{sec:Introduction}

Many applications of industrial interest (e.g., energy and process planning/scheduling and supply‑chain design) often result in large and structured optimization problems. Classical monolithic formulations and solution strategies face computational challenges due to the size and complexity of such problems; as such, practitioners have relied on decomposition algorithms to obtain solutions that are tractable \cite{rahmaniani2017benders, conejo2006decomposition, sahinidis1991convergence}. 

Graph-theoretic modeling has shown to provide a modular framework that facilitates the expression and decomposition of structured optimization problems \cite{cole2025hierarchical, jalving2022graph}. Diverse graph abstractions have been proposed in the literature  \cite{allen2024solution, allen2023multi, allman2019decode, berger2021gboml, berger2021remote, daoutidis2019decomposition, mitrai2020decomposition, mitrai2021efficient, mitrai2024computationally, tang2018optimal, tang2023resolving}; in this work we focus on on the OptiGraph abstraction introduced by Jalving and co-workers and implemented in the open-source package {\tt Plasmo.jl} \cite{jalving2022graph}. Within the OptiGraph, optimization problems are represented as graphs composed of nodes containing optimization subproblems (holding local objectives, constraints, and variables) and edges (encoding linking constraints and variables that couple node subproblems).  

A key observation that inspires our work is that graph modeling provides a flexible and modular foundation for the implementation of decomposition algorithms. Specifically, node subproblems can be isolated and replaced by surrogate models; in addition, tailored solution approaches for the node subproblems can be implemented all without requiring the model code to be rewritten.  The graph approach also provides a framework that enables the implementation of tailored  decomposition approaches \cite{cole2025hierarchical}.

Benders decomposition \cite{benders1962partitioning, rahmaniani2017benders} is a classical and powerful decomposition approach that can be used to tackle problems that have an acyclic graph structure (e.g., trees encountered in stochastic programming). Moreover, we have recently shown that the approach can be applied to problems with general graph structures by conducting graph aggregation operations \cite{cole2025hierarchical}. A key aspect that limits the scalability of Benders decomposition is the repetitive solution of node subproblems.  A promising approach to alleviate this is through the use of surrogate models for node subproblems \cite{bhosekar2018advances, dixit2024decision}. Surrogate models act as effective substitutes for expensive simulations or repetitive  optimization solves. As a result, they can substantially reduce computational effort required in iterative algorithms such as Benders decomposition \cite{semenova2025case}. In our context, surrogates can be used to reduce the computational burden of repeatedly solving graph subproblems. For instance, by learning the behavior of the subproblem solution map/policy, surrogates can quickly estimate optimal values or generate approximate cuts.

Recent work has extensively explored the use of machine learning (ML) models as surrogates to accelerate Benders decomposition across various domains \cite{borozan2024machine, lee2020accelerating, mitrai2024computationally}. These approaches often use neural networks to approximate subproblem solutions and their duals to construct heuristic Benders cuts, employing architectures such standard feed-forward networks \cite{larsen2024fast} or Rectified Linear Unit (ReLU) networks \cite{patel2022neur2sp}. Other approaches leverage more specialized architectures, such as Input Convex Neural Networks (ICNNs) to directly replace the recourse function while maintaining convexity \cite{liu2025icnn}, or Graph Convolutional Networks (GCNs) to generate cuts by exploiting the problem graph structure \cite{neamatian2024neural}. Alternative strategies use reinforcement learning (RL) to generate discrete actions for the master problem \cite{mak2023towards} or apply ML classifiers to enhance cut management strategies within the Benders framework itself \cite{hasan2024accelerating, jia2021benders, mana2023accelerating}.

Despite their success in reducing solution times, ML-based surrogate approaches share limitations that motivate our work. The most significant drawback is the sacrifice of optimality guarantees. These methods are generally approximate/heuristic in nature, providing solutions quickly but forgoing the certificates of optimality that a converged Benders algorithm provides \cite{mak2023towards, larsen2024fast, patel2022neur2sp}. Furthermore, these models often require extensive, offline training of the surrogate model on large datasets \cite{neamatian2024neural, hasan2024accelerating}, and may be constrained by theoretical requirements, such as the convexity needed for ICNNs \cite{liu2025icnn}. 

Multi-parametric programming (mp) can be seen as a surrogate modeling approach that offers a interesting alternative to ML surrogates. By treating external variables and uncertain data as parameters, the entire family of node subproblems can be solved analytically once. The result is a collection of critical regions and associated affine solution functions that can be used to replace subproblem solves with fast critical region look-ups and vector-matrix multiplications to evaluate the affine functions at run time. While the mp approach has been leveraged in other contexts \cite{li2007process, li2008reactive, bemporad2002explicit, avraamidou2020adjustable}, it has yet to be systematically integrated into Benders decomposition of graph-structured problems such as stochastic programs, largely because of scaling issues with an increasing number of scenarios \cite{hugo2005long}. However, recent advances in fast mp solvers \cite{kenefake2022ppopt}, combined with the availability of parallel computing resources and the use of decomposition schemes, enable the use of mp in large and structured problems. For instance, mp has been recently used for enabling the use of model predictive control in a network problems \cite{saini2023noncooperative}.

In this work, we combine graph modeling approach with multi-parametric programming (mp) surrogate models to accelerate  Benders decomposition. We show that the classical Benders optimality cut and the affine value function produced by an mp‑LP are formally identical when evaluated at the same active‑set region. This equivalence guarantees that substituting mp surrogates for exact solves preserves convergence of the Benders algorithm.  We also show that this approach provides insight into algorithmic behavior and enables modular modeling. Specifically, because the mp representation partitions the parameter space into polytopes, it provides an immediate, geometric interpretation of how successive Benders iterations migrate across critical regions, offering a level of transparency and insights that are not attainable with black‑box surrogate models, which can potentially lead to further advancements in decomposition approaches. Moreover, the use of mp surrogates enables the standardization of node subproblems, thus streamlining modeling efforts and enabling plug-and-play of pre-computed mp surrogates.

We implemented our framework using a combination of open-source software tools. We use {\tt Plasmo.jl} to build the optimization model as a graph, {\tt PlasmoBenders.jl} \cite{cole2025hierarchical} to manage the Benders decomposition steps, and a custom Python solver for mp \cite{kenefake2022ppopt,oberdieck2016pop} that connects to the system. We focus our attention on tree structures that appear in stochastic programming problems. The resulting framework allows the modeler to declare stochastic mixed-integer linear programs as graphs, embed mp surrogates into graph nodes. Computational results in a capacity planning problem shows speed‑ups--at times up to 2 orders of magnitude for subproblems--without sacrificing optimality guarantees. To the best of our knowledge, this is the first open-source implementation that integrates graph modeling and mp for Benders decomposition. 

The remainder of the paper is organized as follows. Section \ref{sec:Methods} outlines all the elements of our computational framework. Section \ref{sec:case study} outlines the stochastic programming model used as benchmark. Section \ref{sec:Results and Discussion} discusses computational results and Section \ref{sec:Conclusion} provides concluding remarks and directions for future work.

\section{Methods}
\label{sec:Methods}

\subsection{Graph Modeling}
\begin{figure}
    \centering
    \includegraphics[width=\linewidth]{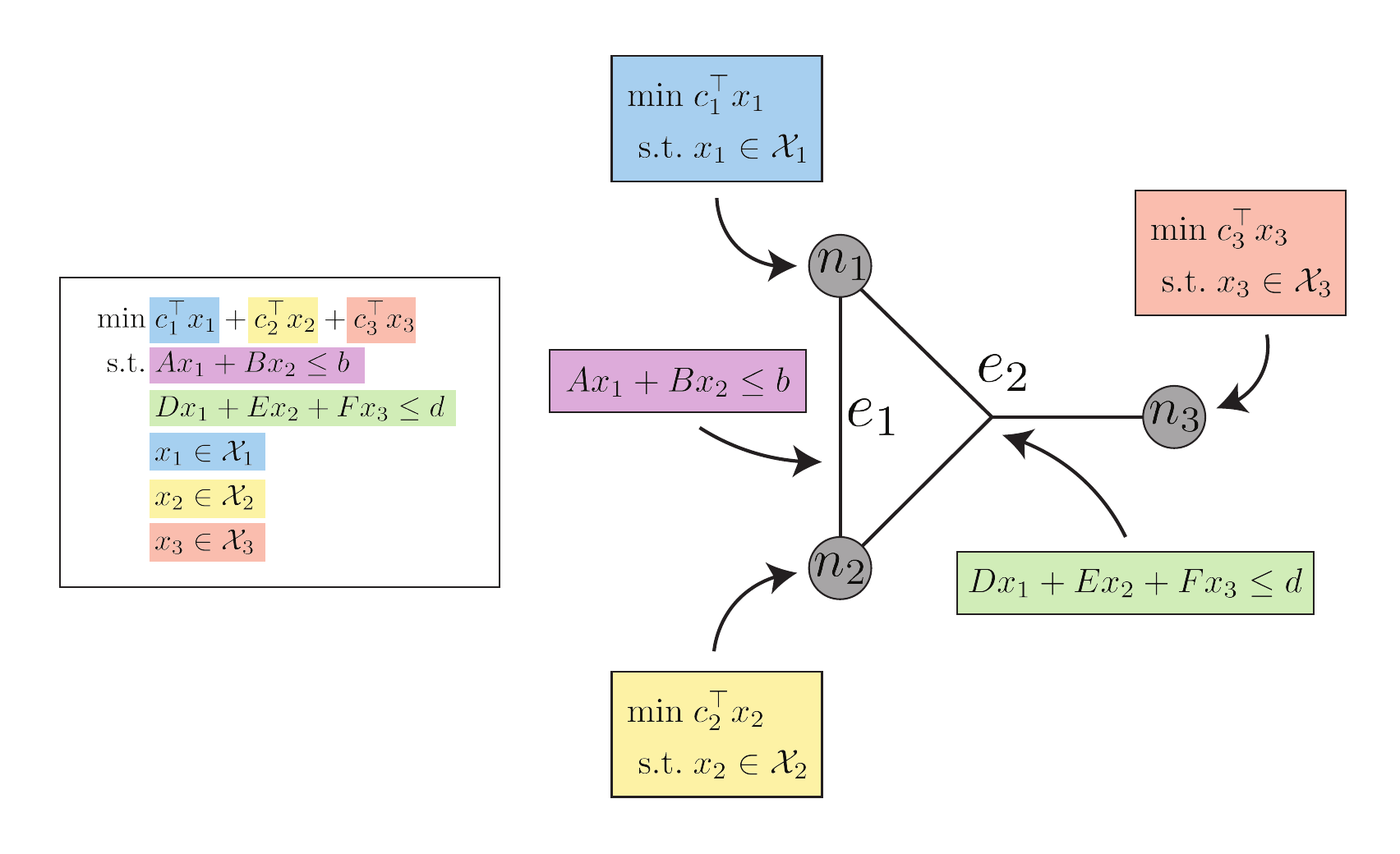}
    \caption{Example of OptiGraph abstraction, where the optimization problem on the left is partitioned into the nodes and edges of the graph on the right.}
    \label{fig:optigraph}
\end{figure}

We represent structure optimization problems as a graph using the OptiGraph abstraction \cite{jalving2022graph} that is implemented in the {\tt Julia} package {\tt Plasmo.jl}. An OptiGraph is made up of a set of OptiNodes, $\mathcal{N}$, and OptiEdges, $\mathcal{E}$. An OptiNode contains an optimization problem (which its own objective, variables, data, and constraints); an OptiEdge contains constraints that link variables accross OptiNodes. The OptiGraph abstraction, $\mg$, can be written as
\begin{subequations}\label{eq:optigraph}
    \begin{align}
        \min_{\{\bx_n\}_{n \in \mn(\mg)}} &\; f\left(\{\bx_n \}_{n \in \mn(\mg)}\right) & (\textrm{Objective}) \label{eq:optigraph_objective} \\
        \textrm{s.t.} &\; \bx_n \in \mathcal{X}_n, \quad n \in \mn(\mg), \quad & (\textrm{Node Constraints}) \label{eq:optigraph_nodes} \\
        &\; g_e(\{\bx_n\}_{n \in \mn(e)})\geq 0, \quad e \in \me(\mg). \quad & (\textrm{Link Constraints}) \label{eq:optigraph_edges}
    \end{align}
\end{subequations}
\noindent where $\mn(\mg)$ are the nodes on $\mg$, $\me(\mg)$ are the edges of $\mg$, $\mn(e)$ are the nodes connecting edge $e$, $\bx_n$ are the decision variables on node $n$, and $\mathcal{X}_n$ is the feasible set for the decision variables on node $n$. An example of the OptiGraph abstraction is shown in Figure \ref{fig:optigraph}. 

The OptiGraph abstraction is able to capture hierarchical relationships within problem structures and form nested subproblems \cite{jalving2022graph, cole2025hierarchical}. The ability to capture hierarchical structures makes the OptiGraph abstraction a flexible modeling tool to embed  surrogate model for node subproblems. In principle, any subproblem node can contain or be replaced by a surrogate model, such as a mp program or a neural network. In this work, we will consider optimization problems where a node is replaced by a mp surrogate. The use of surrogate models to represent node subproblems provides an avenue to standardize the modeling of problems that contain heterogeneous components, as surrogates provide a unifying modeling abstraction for all components. This can help streamline modeling and gain insight into model and algorithmic behavior. 

\subsection{Benders Decomposition}

Benders Decomposition (BD) \cite{benders1962partitioning,rahmaniani2017benders} is a classical, iterative algorithm for solving structured linear programs (LP), mixed-integer linear programs (MILP), and to convex nonlinear problems. In classical BD, an optimization problem with a 2-level tree graph structure is decomposed into a master problem/node and one or more subproblems/subnodes, where the master problem contains a set of complicating variables that, once fixed, enable the parallel solution of subproblems. An iteration of BD includes solving the master problem, passing the primal variable solutions of the master problem to the subproblems, fixing those solutions in the subproblems, solving the subproblems to obtain primal and dual variables that are used for constructing cutting planes that are sent back to the master problem (see Algorithm \ref{alg: benders}). Each iteration of BD results in an upper bound obtained by the solution of the subproblems (a feasible solution) and in a lower bound obtained by  solving the master problem (which contains a relaxation obtained using cutting planes).  

In this work, we apply BD to structured LPs problems of the form \cite{cole2025hierarchical}:
\begin{subequations}\label{eq:BD}
    \begin{align}
        \min &\; c_m^\top x_m + \sum_{w \in W} c_w^\top x_w \\
        \textrm{s.t.} &\; x_m \in \mathcal{X}_m \\
        &\;x_w \in \mathcal{X}_w \quad  w \in W\\
        &\; C_w x_m + D_w x_w \le q_w, \quad w \in W \label{eq:BD_links}.
    \end{align}
\end{subequations}
Here, $x_m$ are the master problem decision variables (complicating/coupling variables), $x_w$ are the subproblem variables for subproblem $w$, $W$ is the set of subproblems, $\mathcal{X}_w$ is the linear feasible region of the subproblem $w$, and ($c_m$, $c_w$, $C_w$, $D_w$, $q_w$) are problem data. The master problem has the form 
\begin{subequations}\label{eq:BD_master}
    \begin{align}
        \underline{\phi}^i_m := \min &\; c_m^\top x_m + \sum_{w \in W} \alpha_w\\
        \textrm{s.t.} &\; x_m \in \mathcal{X}_m \\
        &\; \alpha_w \ge \{ \textrm{cuts} \}, \quad w \in W. \label{eq:NBD_root_cuts}
    \end{align}
\end{subequations}
Here $\alpha_w$ are the subproblem cost variables that provide a lower bound to the optimal solution of the subproblems. $\{\textrm{cuts}\}$ are cutting planes that restrict the value that $\alpha$ can take. Subproblem $w$ takes the form 
\begin{subequations}\label{eq:BD_sub}
    \begin{align}
        \overline{\phi}^i_w(\bar{x}^i_m) := \min &\; c_w^\top x_w\\
        \textrm{s.t.} &\; x_w \in \mathcal{X}_w \label{eq:NBD_subroot_cons} \\
        &\; C_w z_w + D_w x_w \le q_w \label{eq:BD_subproblem_links} \\
        &\; z_w = \bar{x}_m^i \quad (\lambda^i_w). \label{eq:BD_subproblem_dual}
    \end{align}
\end{subequations}
Here, $\bar{x}^i_m$ are the solutions of the master problem at the $i$th iteration, $z_w$ are variables that are fixed to the master problem solution, and $\lambda^i_m$ are dual variables corresponding to constraint \eqref{eq:BD_subproblem_dual}. Note that the constraint \eqref{eq:BD_links} is enforced in the subproblem as \eqref{eq:BD_subproblem_links}. 

Because the master problem solution is enforced on the subproblem, the upper bound at iteration $i$ is given by $UB^i = \min_{j \in \{1, ..., i \}} c^\top_m \bar{x}^j_m + \sum_{w \in W} c_w^\top \bar{x}^j_w$ where $\bar{\cdot}^j$ corresponds to the optimal solution of that variable at iteration $j$. The lower bound at iteration $i$ is given by $LB^i = c^\top_m \bar{x}^i_m + \sum_{w \in W} \bar{\alpha}^i_w$. Finally, Benders cuts are added at each iteration and have the form, for iteration $i$, of $\alpha_w \ge \bar{\phi}_w^i + \lambda^{i \top}_w(x_m - \bar{x}_m^i)$. This approach results in  $|W|$ cuts being added to the master problem at each iteration and is referred to as ``multi-cuts''; an alternative approach is to use a single subproblem cost variable which is restricted by the sum of the right-hand side of the above cut for all subproblems, referred to as an aggregated or single cut.

\begin{algorithm}[H]
\caption{Classical Benders Decomposition}
\label{alg: benders}
\begin{algorithmic}[1]
\State Initialize: iteration counter $i \gets 0$, lower bound $LB \gets -\infty$, upper bound $UB \gets +\infty$
\State Initialize master problem with no Benders cuts
\While{$UB - LB > \epsilon$}
    \State Solve master problem to get $\bar{x}_m^i$
    \For{each scenario/subproblem $w \in W$}
        \State Solve subproblem $w$ with fixed $\bar{x}_m^i$ to get optimal value $\bar{\phi}_w^i$ and duals $\lambda_w^i$
        \State Generate Benders cut: $\alpha_w \ge \bar{\phi}_w^i + \lambda_w^{i\top}(x_m - \bar{x}_m^i)$
        \State Add cut to master problem
    \EndFor
    \State Compute bounds:
    \State \hskip1em $UB \gets \min\{UB, c_m^\top \bar{x}_m^i + \sum_{w \in W} \bar{\phi}_w^i \}$
    \State \hskip1em $LB \gets$ objective value of current master problem
    \State $i \gets i + 1$
\EndWhile
\State \textbf{Return:} Optimal solution $\bar{x}_m$, $\bar{x}_w$ for all $w \in W$
\end{algorithmic}
\end{algorithm}

\subsection{Multiparametric Programming (mp)}
Multiparametric programming (mp) \cite{pappas2021multiparametric,pistikopoulos2020multi} is an  optimization technique that solves optimization problems as explicit functions of input parameters. It has gained significant attention in process systems engineering, particularly in model predictive control \cite{bemporad2002explicit, pistikopoulos2020multi}, multilevel optimization \cite{avraamidou2022multi, avraamidou2019multi}, and robust optimization \cite{avraamidou2020adjustable, lee2023global}. By parameterizing optimization problems, mp enables the derivation of optimal solution maps (for primal and dual variables) as functions of input parameters. Multiparametric programming aims to find the solution of the problem:

\begin{subequations}
\begin{align}
J^{*}(\theta)=&\min_{x\in \mathbb{R}^n}f(x,\theta)\\
\text{s.t.} & ~ g(x,\theta)\leq 0,\\
&~ \theta \in \Theta \subset \mathbb{R}^m
\end{align}
\label{eq: mpproblem}
\end{subequations}

\noindent where $x$ is the vector of decision variables, $\theta$ is the vector of input  parameters, $f(x,\theta)$ is the objective function, $g(x,\theta)$ represents the constraints, and $\Theta$ is a compact polytope representing the bounded uncertain parameters.   The solution map of mp problems are characterized by critical regions (CR), which are subsets of the parameter space where a certain set of constraints is active. When the objective function is linear or quadratic and constraints are linear, these critical regions partition the feasible parameter space into polytopes, and within each region, the optimal solution is described by an affine function of the parameters. This structure allows for efficient implementation of optimal solution strategies, as the solution reduces to a simple function evaluation once the appropriate critical region is identified. 

Here, we assume that the problem is an LP, but it is possible to apply mp to quadratic problems and MILP problems. The optimization formulation of a mp MP problem can be expressed as: 

\begin{subequations}
\begin{align}
z(\theta) = &\min_{x} \ c^\top x + \theta^\top H^\top x \\
\text{s.t.} &~Ax \leq b + F\theta \\
&A_{\text{eq}}x = b_{\text{eq}} + F_{\text{eq}}\theta \\
&\theta \in \Theta := \{ \theta \in \mathbb{R}^q \mid A_{\theta}\theta \leq b_{\theta} \} \\
&x \in \mathbb{R}^n
\end{align}
\label{eq: mpLPproblem}
\end{subequations}


The explicit (mp) solution to eq.~\eqref{eq: mpLPproblem} can be expressed as:

\begin{equation}
 \begin{aligned}
     {x}^* = A^v\theta + b^v ~~ if ~~ CR^v : E^v\theta \leq f^v \\  
     \forall ~~ v = 1, 2, 3, ..., n_{CR}
 \end{aligned}
 \label{eq: mpSol}
 \end{equation}
\noindent where $x^{*}$ is the solution of the problem, $A$ and $b$ represent the affine function defined for each critical region $v$, with $E^v$ and $f^v$ as the corresponding inequality matrices that define the critical region polytopes, and $n_{\text{CR}}$ denotes the total number of critical regions.

\subsection{mp Surrogates for Benders subproblems}
To demonstrate the utility of mp-based surrogates for Benders decomposition, we focus on two-stage stochastic programming  problems (2-level tree graph structure). However, the proposed approach is general and can be applied to Benders decomposition of more general graphs. In a decomposition framework, the problem is split into a master problem and scenario-specific subproblems. For each scenario $k$, the subproblem \eqref{eq:BD_sub} is defined separately and must be solved repeatedly during each iteration of the Benders algorithm. When the subproblem contains scenario-dependent data such as stochastic objective coefficients and right-hand side values, these elements can be treated as parameters. This allows the subproblem to be reformulated as a multiparametric linear program (mp-LP), where the parametric dependence is explicitly modeled. It is worth noting that the same methodology is applicable when the subproblem is a multiparametric mixed-integer linear program (mp-MILP). 

Consider a scenario-specific subproblem \eqref{eq:BD_sub} can be redefined as:

\begin{subequations}\label{eq:BD_sub_scen}
    \begin{align}
        \overline{\phi}^i_{w,k}(\bar{x}^i_m) := \min &\; 
        \begin{bmatrix}
            [c_{w,f}~~\mathbf{0}]\\ 
            [\mathbf{0}~~c_{w,k}] 
        \end{bmatrix}^\top x_w \\
        \text{s.t.} \quad 
        & A_w x_w \leq b_w \label{eq:NBD_subroot_cons_scen} \\
        & C_w z_w + D_w x_w \le 
        \begin{bmatrix}
            q_{w,f}\\
            q_{w,k}
        \end{bmatrix} \label{eq:BD_subproblem_links_scen} \\
        & z_w = \bar{x}_m^i \quad (\lambda^i_w). \label{eq:BD_subproblem_dual_scen}
    \end{align}
\end{subequations}

\noindent Here, $c_{w,f}$ and $q_{w,f}$ are the fixed (non-random) components of the objective and constraints, respectively, while $c_{w,k}$ and $q_{w,k}$ capture the scenario-dependent (random/ uncertain) components. The linking variables $z_w$ are coupled to the master problem solution $\bar{x}_m^i$ through an equality constraint.  This subproblem can be reformulated as a mp-LP where the parameter vector comprises the first stage variables and the uncertain data:
\begin{subequations}\label{eq:BD_sub_scen_mp}
    \begin{align}
        \overline{\phi}^i_{w,k}(\bar{x}^i_m) := \min &\; [c_{w,f}~~\mathbf{0}]^\top x_w + [\mathbf{0}~~c_{w,k}]^\top I^\top x_w \\
        \text{s.t.} \quad 
        & A_w x_w \leq b_w \label{eq:NBD_subroot_cons_scen_mp} \\
        & D_w x_w + C_w z_w \le q_{w,f} + I q_{w,k} \label{eq:BD_subproblem_links_scen_mp} \\
        & z_w = I \bar{x}_m^i \quad (\lambda^i_w). \label{eq:BD_subproblem_dual_scen_mp}
    \end{align}
\end{subequations}

\noindent Here $I$ is the identity matrix. Let $x = \begin{bmatrix}x_w \\ z_w \end{bmatrix}$ be the full decision variable vector, and define the parameter vector as 
\[
\theta = \begin{bmatrix}
    \bar{x}_m^i\\
    c_{w,k}\\
    q_{w,k}
\end{bmatrix}.
\]
Then, the problem in \eqref{eq:BD_sub_scen_mp} fits the standard mp-LP structure described in \eqref{eq: mpLPproblem}, with the following definitions:
\begin{align*}
    &c = \begin{bmatrix} [\mathbf{0}~~c_{w,f}]^\top \\ \mathbf{0} \end{bmatrix}, \quad
    H = 
    \begin{bmatrix}
        \mathbf{0} & \mathbf{0} & \mathbf{0} \\
        I & \mathbf{0} & \mathbf{0}
    \end{bmatrix}, \\
    &A = 
    \begin{bmatrix}
        A_w & \mathbf{0} \\
        D_w & C_w
    \end{bmatrix}, \quad
    b = 
    \begin{bmatrix}
        b_w \\
        q_{w,f}
    \end{bmatrix}, \\
    &F = 
    \begin{bmatrix}
        \mathbf{0} & \mathbf{0} & \mathbf{0} \\
        \mathbf{0} & \mathbf{0} & I
    \end{bmatrix}, \quad
    A_{\text{eq}} = \begin{bmatrix} \mathbf{0} & I \end{bmatrix}, \\
    &b_{\text{eq}} = \mathbf{0}, \quad
    F_{\text{eq}} = \begin{bmatrix} I & \mathbf{0} & \mathbf{0} \end{bmatrix}.
\end{align*}

\noindent
It is important to note that all components used in defining these matrices such as $\mathbf{0}$ and $I$—are themselves matrices of appropriate dimensions. Therefore, the overall matrices \( A, H, F, A_{\text{eq}}, \) and \( F_{\text{eq}} \) are block matrices.

\subsubsection{Equivalence Between Benders Cuts and mp-Based Cuts}

At a given Benders iteration $i$, the master problem \eqref{eq:BD_master} is solved to obtain $\bar{x}_m^i$. The subproblem \eqref{eq:BD_sub} is then solved using $\bar{x}_m^i$ as a fixed input, and the Benders cut is constructed using the subproblem objective value $\bar{\phi}^i$ and dual variable $\lambda^i_w$.

From the mp solution of the subproblem, the parametric objective function is:
\[
\phi(\theta) = c^\top x + \theta^\top H^\top x
\]
For a given parameter vector \( \theta = \bar{\theta} \), which includes the fixed master variable \( \bar{x}_m^i \), the optimal solution lies within a critical region \( v = \bar{v} \), and takes the affine form:
\[
x^* = A^v \bar{\theta} + b^v
\]

Substituting this into the objective function gives:
\begin{equation}\label{eq:phi_mp}
\bar{\phi}^i = c^\top (A^v \bar{\theta} + b^v) + \bar{\theta}^\top H^\top (A^v \bar{\theta} + b^v)
\end{equation}

We now compute the gradient of the value function with respect to the master decision variable \( x_m \), embedded in \( \theta \). Denote by \( \theta = \begin{bmatrix} \bar{x}_m \\ \theta^* \end{bmatrix} \), where \( \theta^* = \begin{bmatrix} c_{w,k} \\ q_{w,k} \end{bmatrix} \), and recall that only the first block \( \bar{x}_m \) is a variable in the master problem as \( x_m \).
We compute:
\[
\nabla_{x_m} \phi(\theta) = \nabla_{\theta} \phi(\theta) \cdot \nabla_{x_m} \theta
\]

But since \( \theta = [\bar{x}_m; \theta^*] \), only the top block of \( \theta \) is affected, so:
\[
\nabla_{x_m} \phi(\theta) = (c + \theta^\top \bar{H}^\top) A^v_{x_m}
\]
where \( A^v_{x_m} \) is the portion of \( A^v \) multiplied by \( x_m \) and $\bar{H} = \begin{bmatrix} I & \mathbf{0} \\ \mathbf{0} & \mathbf{0} \end{bmatrix}$.

This leads to the mp-based affine cuts:
\begin{equation}\label{eq:Benders_cut_mp}
\theta_w \geq \bar{\phi}^i + \left( c + \theta^{*\top} \bar{H}^\top \right) A^v_{x_m} (x_m - \bar{x}_m^i)
\end{equation}

Now recall the standard Benders cut derived from duality:
\[
\theta_w \geq \bar{\phi}^i + \lambda^{i \top}_w (x_m - \bar{x}_m^i)
\]

We now formally establish that both cuts are equivalent under the framework of mp using active-set theory. We aim to show that:
\[
\lambda_w^i = \left( c + \theta^{*\top} \bar{H}^\top \right) A^v_{x_m}
\] 

For a given critical region, the solution is determined by the active set of constraints. Revisiting the subproblem \eqref{eq:BD_sub_scen_mp}, and its mp-LP representation \eqref{eq: mpLPproblem}, the solution within a critical region can be equivalently expressed as the following active-set-based mp-LP:

\begin{subequations}
\begin{align}
z_{AS}(\theta) = &\min_{x} \ c^\top x + \theta^\top H^\top x \\
\text{s.t.} \quad &A_{AS}x = b_{AS} + F_{AS}\theta \\
&\theta \in \Theta := \{ \theta \in \mathbb{R}^q \mid CR^A \theta \leq CR^b \}, \quad x \in \mathbb{R}^n
\end{align}
\label{eq: mpLPproblem_active_set}
\end{subequations}

Here, $A_{AS} = \begin{bmatrix} A_k \\ A_{\text{eq}} \end{bmatrix}$, $b_{AS} = \begin{bmatrix} b_k \\ b_{\text{eq}} \end{bmatrix}$, and $F_{AS} = \begin{bmatrix} F_k \\ F_{\text{eq}} \end{bmatrix}$, where $A_k$, $b_k$, and $F_k$ are the rows of the original matrices active at optimality.

The Lagrangian for this equality-constrained mp-LP is:
\begin{equation}
\begin{aligned}
\nabla_x \mathcal{L}(x, \lambda, \theta) 
&= \nabla_x \left( (H\theta + c)^\top x + \lambda^\top (A_{AS} x - b_{AS} - F_{AS} \theta) \right) \\
&= H\theta + c + A_{AS}^\top \lambda
\end{aligned}
\label{eq: Lagrangian_function}
\end{equation}

The Karush-Kuhn-Tucker (KKT) conditions for this region are:
\begin{equation}
\begin{aligned}
& H\theta + c + A_{AS}^\top \lambda = 0 \\
& A_{AS}x = b_{AS} + F_{AS}\theta.
\end{aligned}
\end{equation}

From the KKT conditions of the active-set mp-LP, we get:
\[
H\theta + c + A_{AS}^\top \lambda = 0 \quad \Rightarrow \quad \lambda = - (A_{AS}^\top)^{-1} (H\theta + c)
\].

The resulting affine primal solution is:
\[
x = A_{AS}^{-1}(F_{AS} \theta + b_{AS}) = A^v \theta + b^v
\].

The dual of \eqref{eq: mpLPproblem_active_set} is:
\begin{equation}
\begin{aligned}
z_{AS}(\theta) = &\max_{\lambda} \ \lambda^\top(b_{AS} + F_{AS}\theta) \\
\text{s.t.} \quad &A_{AS}^\top \lambda \geq c + H\theta, \quad \lambda \geq 0, \quad \theta \in \Theta.
\end{aligned}
\label{eq: mpLPproblem_active_set_dual}
\end{equation}

From the dual expression we get the following:
\[
\lambda^\top (b_{AS} + F_{AS} \theta) = - (H\theta + c)^\top (A_{AS}^\top)^{-1} (b_{AS} + F_{AS} \theta)
\]
\[
= (c + \theta^\top H^\top)^\top A^v \theta + \text{const.}
\]

Now, focusing only on how the objective varies with \( x_m \) (i.e., the dual variable of \( x_m \)) we see:
\[
\lambda^\top (F_{AS})_{x_m} = (c + \theta^\top H^\top) A^v_{x_m}
\]

Therefore:
\[
\lambda^{i \top}_w = \left( c + \theta^{*\top} \bar{H}^\top \right) A^v_{x_m}
\]

which completes the proof that the cuts are identical:
\[
\bar{\phi}^i + \lambda^{i\top}_w (x_m - \bar{x}_m^i)
=
\bar{\phi}^i + \left( c + \theta^{*\top} \bar{H}^\top \right) A^v_{x_m} (x_m - \bar{x}_m^i)
\]
\subsection{Software Framework}

The algorithmic framework proposed was implemented in open source software available at \url{https://github.com/Avraamidou-Research-Group-CESE/MP-Plasmo-Benders}. The framework relies on Python (for the mp-surrogates) and Julia (for graph-based modeling and Benders decomposition). The integration of two different programming environments highlights the modularity that is enabled by using mp surrogates within a graph modeling environment. Below we outline the implementation of the framework. 

We used the PPOPT toolbox \cite{kenefake2022ppopt} in Python to generate the mp solution. The subproblem \ref{eq:BD_sub_scen_mp} can be formulated and solved as a mp problem after defining parameters. The solution of the resulting mp problem contains $A^v$ and $b^v$ that represent the affine function defined for each critical region $v$, with $E^v$ and $f^v$ as the corresponding inequality matrices that define the critical region as shown in eq. \ref{eq: mpSol}. This list of matrices for each critical region, representing the entire map of optimal solutions of the subproblem \ref{eq:BD_sub_scen_mp} for any realization of the defined parameters, is saved and loaded to Julia.

\begin{figure}[H]
    \centering
    \includegraphics[width=0.7\linewidth]{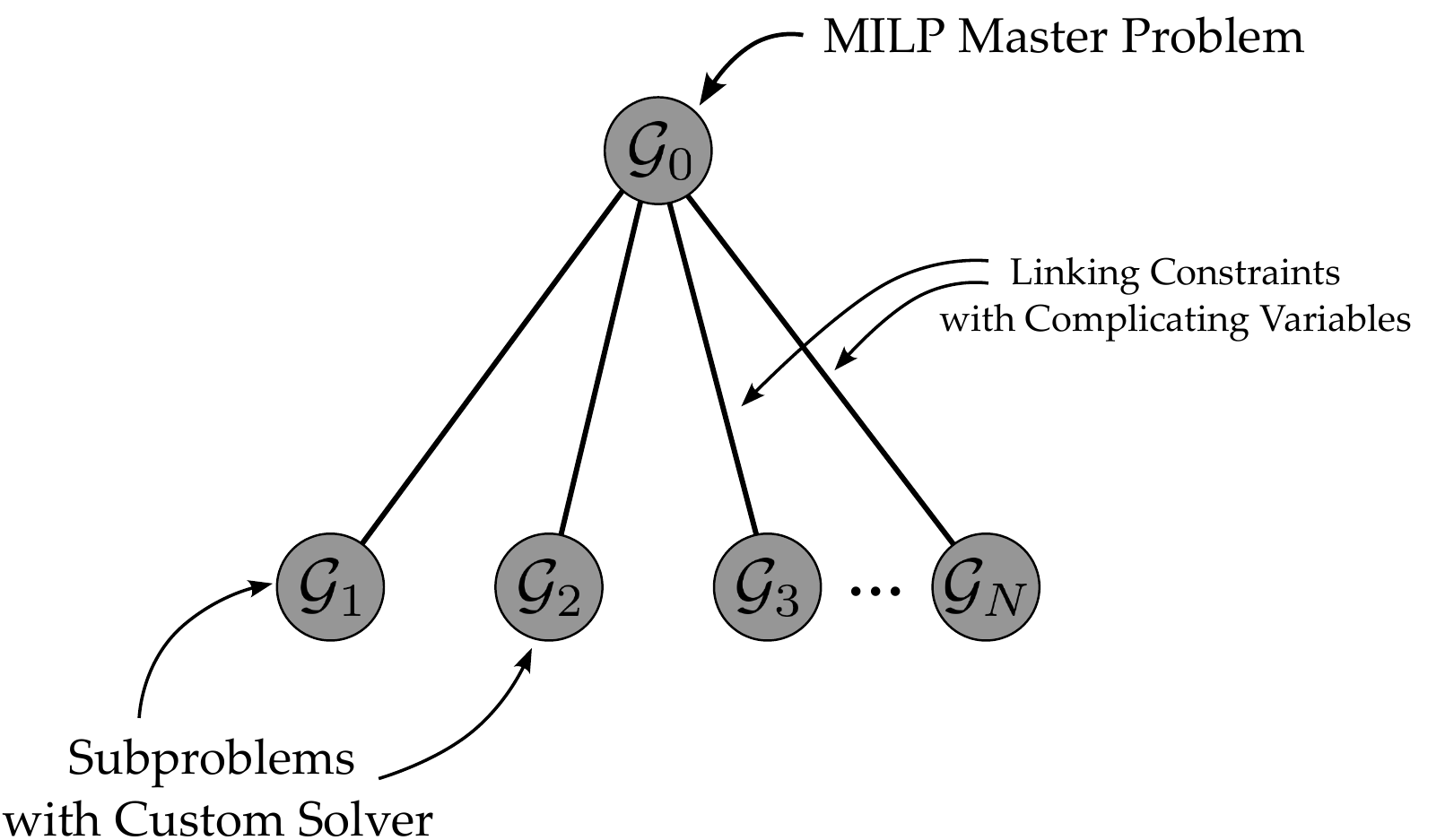}
    \caption{Graph representation of the problem solved with Benders. The master problem is stored on a subgraph, $\mathcal{G_0}$, and each subproblem is stored on  independent subgraphs $\mathcal{G}_1$ - $\mathcal{G}_N$. Edges contain the constraints with complicating variables. }
    \label{fig:graph_representation}
\end{figure}

After loading the mp model solution, we build the full optimization problem in {\tt Plasmo.jl} \cite{jalving2022graph}. The master problem and each subproblem are placed on independent subgraphs with linking constraints between the master problem and each subproblem (Figure \ref{fig:graph_representation}). Each subproblem is independent of the other subproblems such that each linking constraint only links the master problem subgraph and each subproblem subgraph. Once this structure has been created, Benders decomposition can be applied via the package {\tt PlasmoBenders.jl}. This package takes the user-defined OptiGraph structure and requires the user to set the master subgraph. Once this subgraph is set, {\tt PlasmoBenders.jl} detects the downstream subgraphs, fixes the complicating variables in the subproblem subgraphs, and adds the cost to go variable(s). In the case of Figure \ref{fig:graph_representation}, the user can set subgraph $\mathcal{G}_0$ as the master problem and {\tt PlasmoBenders.jl} will recognize that subgraphs $\mathcal{G}_1$ - $\mathcal{G}_N$ are the subproblem subgraphs.

To implement the mp surrogate models within the graph modeling framework, we define a custom solver in Julia such that no changes are required to the existing software of {\tt Plasmo.jl} and~{\tt PlasmoBenders.jl}. This is a ``custom solver'' in the sense that we have defined our own methods in the Julia code that uses the mp data to solve the optimization problem and query primal-dual information, independent of typical LP solvers such as HiGHS or Gurobi. This custom solver interfaces with {\tt JuMP.jl}'s~{\tt MathOptInterface.jl} (used in the backend of {\tt Plasmo.jl}) to query primal-dual information from mp surrogate data that was determined from PPOPT toolbox. The user configures the solver as the optimizer for the subproblem corresponding to a specific subgraph and subsequently retrieves the mp solution data directly from the solver. Once {\tt optimize!} is called on the subproblem, primal and dual information queried from the solver will be returned internally from the mp solution functions. By creating the custom solver, we allow for flexibly working with existing, registered packages and provide a framework for handling surrogate modeling approaches within the graph-based modeling framework. 

The overall outline of the proposed framework is summarized in algorithm \ref{alg: mp-Plsmo-Benders} below.

\begin{algorithm}[H]
\caption{Benders Decomposition with mp Surrogate via \texttt{PlasmoBenders.jl}}
\label{alg: mp-Plsmo-Benders}
\begin{algorithmic}[1]

\State \textbf{mp Surrogate Generation}
  \begin{itemize}
    \item Formulate the mp subproblem \eqref{eq:BD_sub_scen_mp}.
    \item Solve it using the PPOPT toolbox to obtain critical regions and associated affine solutions.
    \item Extract and store the surrogate model data: $(E^v, f^v, A^v, b^v)$.
  \end{itemize}

\State \textbf{Initialization in \texttt{PlasmoBenders.jl}}
  \begin{itemize}
    \item Formulate the master problem and create a custom Benders solver.
    \item Load the mp surrogate model data into the custom solver.
    \item Initialize model parameters for the optimization problem.
  \end{itemize}

\State \textbf{Run \texttt{PlasmoBenders.jl} Solver with mp Surrogate}
  \begin{itemize}
    \item Execute the solver. Internally, it performs an iterative Benders decomposition loop using the mp surrogate model.
  \end{itemize}

  \Statex \hspace{1em} The internal loop proceeds as follows:

  \While{$UB^i - LB^i > \epsilon$ \textbf{(handled internally)}}
    \State Solve the master problem~\eqref{eq:BD_master} to obtain current solution $\bar{x}_m^i$.
    \For{each subproblem node $k$ and scenario $w$}
      \State Construct parameter vector $\bar{\theta} = (\bar{x}_m^i, c_{w,k}, q_{w,k})$.
      \State Locate the critical region corresponding to $\bar{\theta}$.
      \State Evaluate subproblem cost $\bar{\phi}^i$ using affine surrogate \eqref{eq:phi_mp}.
      \State Add Benders cut \eqref{eq:Benders_cut_mp} to the master problem.
    \EndFor
    \State Update bounds $LB^i$ and $UB^i$.
  \EndWhile

\end{algorithmic}
\end{algorithm}

Replacing LP solves with simple affine evaluations makes the
sub‑problem phase almost cost‑free while preserving convergence
guarantees, as the cuts are unchanged.


\section{Case Study: Capacity–Expansion Planning Under Uncertainty}
\label{sec:case study}

To highlight how the proposed methodology and software framework can be applied, we present a case study of a capacity-expansion planning problem under uncertainty for a chemical process \cite{sahinidis1989optimization}. The capacity of the chemical process units can be expanded at different points in time, while there is market uncertainty in the expected demands, availabilities, and prices of chemicals. The first-stage (master problem) variables represent capacity expansion decisions. The subproblem variables capture the purchase of raw materials, as well as the production and sale of intermediate and final products. While this case study is specific to a chemical process, the mathematical structure of this problem is common in many applications and domains, including power infrastructure planning \cite{lara2018deterministic} and energy systems planning \cite{jacobson2024computationally}.

\subsection{Problem Description}
In its most basic form, the long-term capacity expansion and production planning problem with fixed-charge expansion costs can be stated as follows. Given an extended network of continuous production processes, \( p \), and chemicals, \( j \), the primary goal is to determine the optimum selection of processes and how their capacities should be expanded over a future horizon consisting of a set of time periods, \( t \). Over this planning horizon, the demands, availabilities and prices of the chemicals are forecast as stochastic quantities, while the process capital investment, operating costs and expansion bounds, are deterministically known. Information is also given on how each process converts the raw materials into products in terms of conversion factors.

Because of the presence of uncertainty in the market conditions, the expected demands, availabilities and prices of the chemicals are predicted according to a set of scenarios, \( k \in \mathcal{K} \). To manage this uncertainty, the decisions are grouped into a couple of  stages. In the absence of knowing the uncertainty, the capacity expansion decisions have to be made in the \emph{first stage}, irrespective of which scenario occurs. However, the production rates of the processes and the amounts of chemicals purchased and sold can be manipulated in the \emph{second stage} to best satisfy the market conditions revealed under each scenario. 

The objective is to identify the optimal capacity plan to implement that maximizes the expected profit while still ensuring feasible plant operations over all the uncertain future realizations. The decision whether an expansion should take place or not is modeled through the use of a binary variable. The model is, therefore, formulated as a MILP problem (CEP) with the notation presented in eq. \eqref{eq:orig}.

\begin{subequations} \label{eq:orig}
\begin{align}
\min_{x_{pt},y_{pt},q_{pt}}
    & \;
      \sum_{p\in\mathcal{P}}\sum_{t\in\mathcal{T}}\left[ \alpha_{pt}\,x_{pt} + \beta_{pt} y_{pt} \right]
      \;+\;
      \mathbb{E}_{\xi}\!\bigl[V(q,\xi)\bigr]
      \label{eq:orig_obj} \\[4pt]
\text{s.t. }\;
    & E_{pt}^{L} y_{pt} \;\le\;
      x_{pt}
      \;\le\; E_{pt}^{U} y_{pt}
      &&\forall p\in\mathcal{P},\;t\in\mathcal{T}
      \label{eq:orig_cap} \\[2pt]
    & 0 \;\le\; q_{pt} \;\le\; Q_{pt}^{U}
      &&\forall p\in\mathcal{P},\;t\in\mathcal{T}
      \label{eq:orig_bounds} \\
    & q_{pt} - q_{pt-1} - x_{pt} = 0 &&\forall p\in\mathcal{P},\;t\in\mathcal{T} \label{eq:orig_cuml} \\
    & x_{pt} \geq 0 &&\forall p\in\mathcal{P},\;t\in\mathcal{T} \label{eq: orig_Non_neg}\\
    & y_{pt} \in \{0,1\} &&\forall p\in\mathcal{P},\;t\in\mathcal{T} \label{eq: orig_Binary_vars}
\end{align}
\end{subequations}

\noindent Here, $x_{pt}$ are capacity expansion decisions, $q_{pt}$ are cumulative expansion decisions, and $y_{pt}$ are binary decisions to expand or not at given time. The objective \eqref{eq:orig_obj} is the minimization of overall cost, considering also the expected recourse cost. Constrains \eqref{eq:orig_cap} restricts within physical limits the capacity expansion of processes. In the same constrains, the binary variables, $y_{pt}$, denote the time $t$ that capacity expansion takes place. Constrains \eqref{eq:orig_bounds} bound the total available capacity to the maximum amount allowed ($Q^U_{pt}$) and enforce nonnegative capacities. Equality constrains \eqref{eq:orig_cuml} add new expansions to the already available capacities from pervious time periods. Finally, constraints \ref{eq: orig_Non_neg} enforce non-negativity for capacity variables, and constraints \ref{eq: orig_Binary_vars} define binary indicator constraints.

For a fixed first‑stage plan \(q_{pt}\) and realization
\(\xi=(A_{jt}^{k},D_{jt}^{k},\gamma_{jt}^{k},\varphi_{jt}^{k})\), the recourse cost variable, $V(q,\xi)$, is calculated based on the  optimization problem \eqref{eq:V_def}.
\begin{subequations}\label{eq:V_def}
\begin{align}
V(q,\xi)\;=\;
\min_{\substack{p_{pt},\,b_{jt},\,s_{jt}}}
\;&
  -\sum_{t\in\mathcal{T}} 
  \left[
      -\sum_{p\in\mathcal{P}}\sigma_{pt}\,p_{pt}
    + \sum_{j\in\mathcal{J}}\gamma_{jt}\,b_{jt}
    - \sum_{j\in\mathcal{J}}\varphi_{jt}\,s_{jt}
  \right]
\label{eq:V_def_obj}
\\
\text{s.t. }\;
& p_{pt} - q_{pt} \leq 0 &&\forall p,t \label{eq:V_def_c1} \\
& s_{jt} \leq A_{jt},\quad b_{jt} \leq D_{jt} &&\forall j,t \label{eq:V_def_c2} \\
& p_{pt} \geq 0,\quad s_{jt}, b_{jt} \geq 0 &&\forall p,j,t \label{eq:V_def_c3} \\
&
  \sum_{p\in\mathcal{P}}(\mu_{jp}-\eta_{jp})p_{pt}
  + b_{jt}-s_{jt} \;=\; 0
  &&\forall j,t \label{eq:V_def_c4}
\end{align}
\end{subequations}
Once market and technological uncertainties are realized, the plant operates at recourse production levels \(p_{pt}^{k}\), with external purchases \(b_{jt}^{k}\) and sales \(s_{jt}^{k}\) of intermediates and products. These decisions are determined as the optimal solution of the subproblem~\eqref{eq:V_def}. The recourse function is defined by the objective~\eqref{eq:V_def_obj}, which minimizes the negative of the profit. Constraints~\eqref{eq:V_def_c1} ensure that production does not exceed the available capacity. Constraint~\eqref{eq:V_def_c2} restricts the quantities of intermediates and final products that can be bought or sold. Finally, constraint~\eqref{eq:V_def_c4} enforces material balance by ensuring that the total amount produced and purchased equals the amount consumed or sold.

\subsection{Benders Reformulation}
\label{ssec:BD_multicut}

To embed the stochastic capacity-expansion planning framework with master problem \eqref{eq:orig} and recourse model \eqref{eq:V_def} into a Benders decomposition framework, a couple of variants are considered: a \emph{multi-cut} scheme, which generates one cut for each scenario–period pair \((k,t)\in\mathcal{K}\times\mathcal{T}\) per iteration, and a \emph{single-cut} scheme, which aggregates all cuts into a single inequality. Let \(\pi_k\) denote the probability of scenario \(k\), and assume uniform probability \(\pi_k = \frac{1}{|\mathcal{K}|}\) for all \(k\in\mathcal{K}\). The formulation is designed to accommodate a relatively complete recourse structure for this problem. Below the problem formulation of the master and subproblems is presented.

Master problem at iteration $l$ can be given by formulation \ref{eq:master_cuts_multicut}. For BD, we introduce auxiliary variables \(\alpha_{kt}\) (multi-cut) or a single variable \(\alpha\) (single-cut) to approximate the recourse cost. The master problem in the multi-cut form is then given by:

\begin{subequations}\label{eq:master}
\begin{align}
\min_{x,y,q,\alpha_{kt}} & \quad
    \sum_{p\in\mathcal{P}}\sum_{t\in\mathcal{T}} \bigl[\alpha_{pt}\,x_{pt} + \beta_{pt}\,y_{pt}\bigr]
  + \sum_{k\in\mathcal{K}} \sum_{t\in\mathcal{T}} \pi_k\,\alpha_{kt}
  \label{eq:master_obj_multicut} \\[4pt]
\text{s.t.} \quad
& E_{pt}^{L} y_{pt} \le x_{pt} \le E_{pt}^{U} y_{pt}
    && \forall p,t \label{eq:master_cap} \\
& 0 \le q_{pt} \le Q_{pt}^{U}
    && \forall p,t \label{eq:master_bounds} \\
& q_{pt} - q_{p,t-1} - x_{pt} = 0
    && \forall p,t \label{eq:master_link} \\
& y_{pt} \in \{0,1\}, \quad x_{pt} \ge 0
    && \forall p,t \label{eq:master_domain} \\
& \alpha_{kt} \ge \bar{\phi}^{\,\ell}_{kt}
    + {\lambda_{kt}^{\,\ell}}^{\!\top} (q_{pt} - \bar{q}_{pt}^{\,\ell})
    && \forall k,t,\; \ell = 1,\dots,i \label{eq:master_cuts_multicut}
\end{align}
\end{subequations}

In the single-cut formulation, only the objective and the cut constraint are modified as follows:\\
Objective (replace \eqref{eq:master_obj_multicut}):
  \[
  \min_{x,y,q,\alpha} \quad
    \sum_{p\in\mathcal{P}}\sum_{t\in\mathcal{T}} \bigl[\alpha_{pt}\,x_{pt} + \beta_{pt}\,y_{pt}\bigr]
  + \alpha
  \]\\
Cut constraint (replace \eqref{eq:master_cuts_multicut}):
  \[
  \alpha \ge \sum_{k\in\mathcal{K}} \sum_{t\in\mathcal{T}} \pi_k
  \left[
    \bar{\phi}^{\,\ell}_{kt} + {\lambda_{kt}^{\,\ell}}^{\!\top}(q_{pt} - \bar{q}_{pt}^{\,\ell})
  \right]
  \quad \forall \ell = 1,\dots,i
  \]

Here \(\bar{q}_{pt}^{\,\ell}\) is master solution at iteration \(\ell\), \(\bar{\phi}_{kt}^{\,\ell}\) is optimal value of the \((k,t)\) subproblem at iteration \(\ell\), and \(\lambda_{kt}^{\,\ell}\) is dual multipliers from the subproblem at iteration \(\ell\).
Single‑period sub‑problem (scenario $k$, period $t$) is given by formulation \ref{eq:subproblem_kt}.
Given the master solution \(\bar{q}_{pt}^{\,\ell}\), the recourse LP for
each scenario–period pair \((k,t)\) is

\begin{subequations}\label{eq:subproblem_kt}
\begin{align}
\bar{\phi}^{\,\ell}_{kt}:=\;
\min_{\substack{p_{p},\,b_{j},\,s_{j}}}\;
  & -\sum_{p\in\mathcal{P}}\sigma_{pt}\,p_{p}
    + \sum_{j\in\mathcal{J}}\gamma_{jt}^{k}\,b_{j}
    - \sum_{j\in\mathcal{J}}\varphi_{jt}^{k}\,s_{j} 
                                  \\
\text{s.t. }&
    p_{p} \;\le\; \bar{q}_{pt}^{\,\ell}
    &&\forall p \;\;                      \label{eq:sub_cap}\\
&
    s_{j} \;\le\; A_{jt}^{k},\quad
    b_{j} \;\le\; D_{jt}^{k}
    &&\forall j \;\;                  \label{eq:sub_bounds}\\
&
    \sum_{p}(\mu_{jp}-\eta_{jp})\,p_{p}
      + b_{j}-s_{j}=0
    &&\forall j \;\;               \label{eq:sub_balance}\\
&
    p_{p},\,b_{j},\,s_{j}\ge 0
    &&\forall p,j .
\end{align}
\end{subequations}

\subparagraph{mp reformulation of the sub‑problem.}

Problem~\eqref{eq:subproblem_kt} can be written in the compact
multi‑parametric LP form
\[
\min_{x}\; c^{\!\top}x + \theta^{\!\top} H^{\!\top}x
\quad\text{s.t.}\quad
Ax \le b + F\theta,\;\;
x\ge 0,
\]
by defining
\[
x :=
\begin{bmatrix}
p \\[2pt] b \\[2pt] s
\end{bmatrix},
\qquad
\theta :=
\begin{bmatrix}
\bar{q}_{p t}^{\,\ell}\\[2pt]
\gamma_{j t}\\[2pt]
\sigma_{p t} \\[2pt]
\varphi_{j t}\\[2pt]
A_{j t}\\[2pt]
D_{j t}^{k}
\end{bmatrix}.
\]
The mp solver partitions
the parameter space \(\Theta=\{\theta\mid k\in\mathcal{K}\}\) into
critical regions; inside each region the optimal decision and cost are
piecewise affine:
\[
x^{*}(\theta)\;=\;A^{v}\theta + b^{v},
\qquad
\bar{\phi}_{kt}(\theta)\;=\;\alpha^{v}{}^{\!\top}\theta + \beta^{v}.
\]

\subsection{Illustrative Example and Numerical Experiments Setup}

\begin{figure}[H]
   \centering
 \includegraphics[width=.6\linewidth]{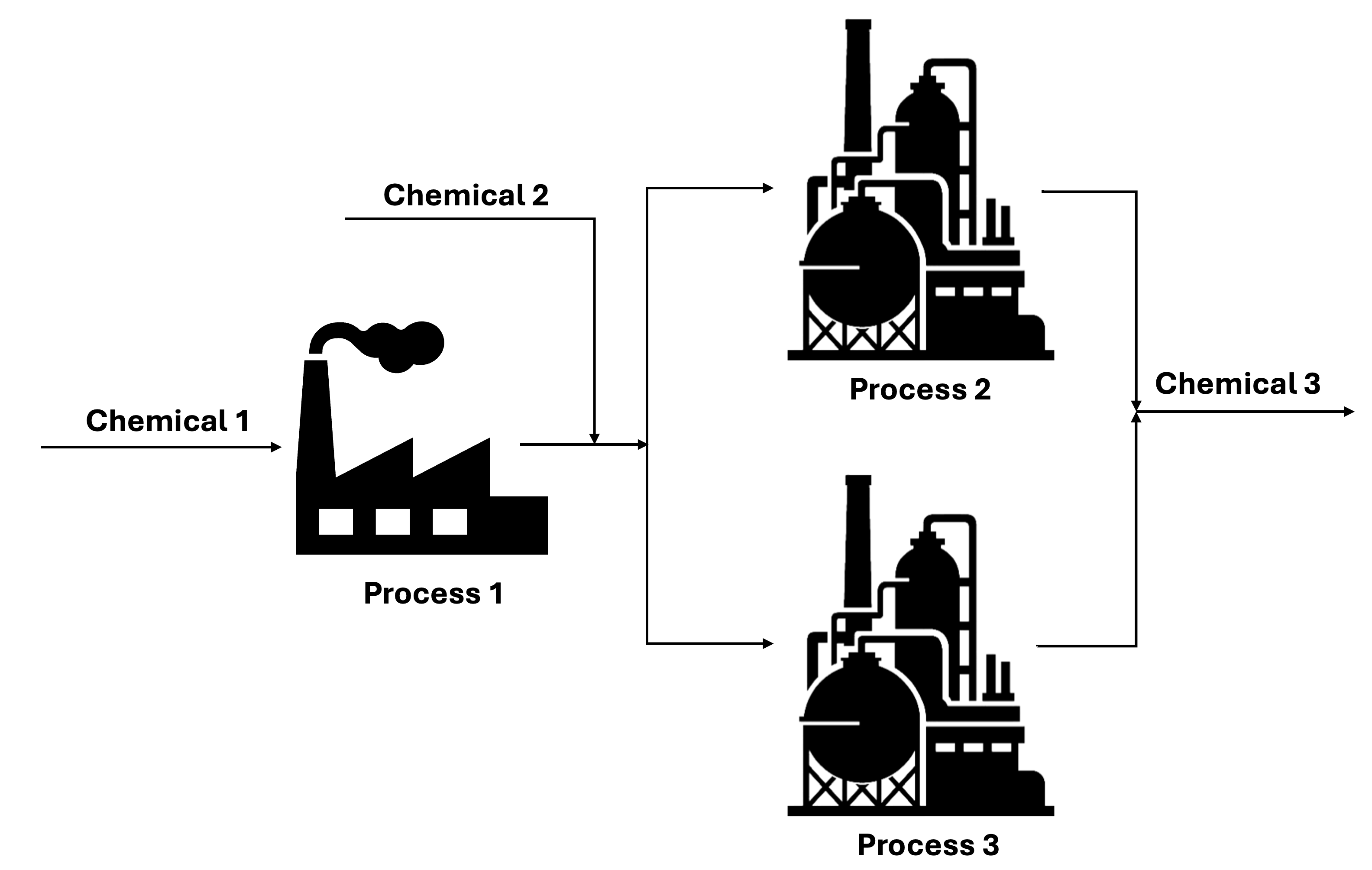}
   \caption{Superstructure of the plant transforming
            \textit{Chemical 1} and \textit{Chemical 2} into the final
            product \textit{Chemical 3}. \cite{iyer1998bilevel}}
   \label{fig:superstructure}
\end{figure}

The network used to illustrate the stochastic capacity expansion problem is adapted from \cite{iyer1998bilevel} and shown in Fig.~\ref{fig:superstructure}. It consists of three continuous production processes and three chemicals operated over a finite planning horizon of \(T\) discrete time periods (\(t = 1,\dots,T\)). For the numerical experiments, we consider four planning horizon lengths: \(T = 5, 10, 25, 50\), and five sizes for the number of scenarios: \(K = 50, 200, 500, 1000, 2000\).

The stochasticity in the model arises from uncertain market and technological conditions—specifically, the demands, availabilities, and prices of the chemicals. Each of these three chemical-related quantities is considered uncertain, yielding a total of six uncertain parameters (i.e., \(3 \times 2 = 6\)). In addition, the operating cost parameters \(\sigma_{pt}\) for each of the three processes and the first-stage cumulative capacity parameters \(q_{pt}\) are also treated as inputs to the model. This brings the total number of parameters used in the multiparametric (mp) formulation of the subproblem~\eqref{eq:subproblem_kt} to twelve.

Although in practice prices, demands, and availabilities are often statistically correlated, we assume for simplicity that they are mutually independent. Each scenario \(k\in\mathcal{K}\) is treated as an independent realization of uncertainty with a corresponding probability weight \(\pi_k\), which is used in the expectation operator in the stochastic model.

For completeness, Table~\ref{tab:mean_values_t4} lists the mean values employed at period \(t=4\), where the rows correspond to processes \(p=1,2,3\) and the columns to time periods \(t=1,\dots,4\). To generate the full scenario set, we sample values from normal distributions centered at these means, using 10\% of the mean as the standard deviation for each distribution.\\

\begin{table}[H]
\footnotesize
\centering
\caption{Mean parameter values used for the base case ($p=1,2,3$; $t=1\ldots4$)}
\begin{tabular}{lcccc}
\toprule
           & $t=1$ & $t=2$ & $t=3$ & $t=4$ \\
\midrule
$\alpha_{pt}$ & 1.38 & 1.67  & 2.22  & 3.58  \\
              & 2.72 & 3.291 & 4.381 & 7.055 \\
              & 1.76 & 2.13  & 2.834 & 4.565 \\
\midrule
$\beta_{pt}$  & 85   & 102.85 & 136.89 & 220.46 \\
              & 73   & 88.33  & 117.56 & 189.34 \\
              & 110  & 133.10 & 177.15 & 285.31 \\
\midrule
$\sigma_{pt}$ & 0.40 & 0.48 & 0.64 & 1.03 \\
              & 0.60 & 0.72 & 0.96 & 1.55 \\
              & 0.50 & 0.60 & 0.80 & 1.29 \\
\midrule
$A_{jt}^{k}$  & 6.00  & 7.26  & 9.66  & 15.56 \\
              & 20.00 & 24.20 & 32.21 & 51.87 \\
              & 0.00  & 0.00  & 0.00  & 0.00  \\
\midrule
$D_{jt}^{k}$  & 0.00  & 0.00  & 0.00  & 0.00  \\
              & 0.00  & 0.00  & 0.00  & 0.00  \\
              & 30.00 & 36.30 & 48.31 & 77.81 \\
\midrule
$\gamma_{jt}^{k}$ & 0.00  & 0.00  & 0.00  & 0.00  \\
                  & 0.00  & 0.00  & 0.00  & 0.00  \\
                  & 26.20 & 31.70 & 42.19 & 67.95 \\
\midrule
$\varphi_{jt}^{k}$   & 4.00 & 4.84  & 6.44  & 10.37 \\
                  & 9.60 & 11.61 & 15.46 & 24.90 \\
                  & 0.00 & 0.00  & 0.00  & 0.00  \\
\bottomrule
\end{tabular}
\label{tab:mean_values_t4}
\end{table}

\vspace{1ex}
\noindent\textbf{Deterministic bounds and stoichiometry}

\begin{table}[H]
\footnotesize
\centering
\begin{tabular}{lccc}
\toprule
               & $p=1$ & $p=2$ & $p=3$ \\
\midrule
$E_{pt}^{L}$   & 1  & 10 & 10 \\
$E_{pt}^{U}$   & 6  & 30 & 30 \\
$Q_{pt}^{U}$   & 100 & 100 & 100 \\
\bottomrule
\end{tabular}
\caption{Lower and upper cumulative capacity bounds ($E_{pt}^{L},E_{pt}^{U}$) and
         per‑period expansion limit ($Q_{pt}^{U}$).}
\label{tab:bounds}

\end{table}

\begin{table}[H]
\footnotesize
\centering
\begin{tabular}{lccc}
\toprule
               & $p=1$ & $p=2$ & $p=3$ \\
\midrule
\multicolumn{4}{l}{\textit{Consumption coefficients }$\eta_{jp}$} \\[2pt]
$j=1$          & 1.11 & 0    & 0    \\
$j=2$          & 0    & 1.22 & 1.05 \\
$j=3$          & 0    & 0    & 0    \\
\midrule
\multicolumn{4}{l}{\textit{Yield coefficients }$\mu_{jp}$} \\[2pt]
$j=1$          & 0 & 0 & 0 \\
$j=2$          & 1 & 0 & 0 \\
$j=3$          & 0 & 1 & 1 \\
\bottomrule
\end{tabular}
\caption{Stoichiometric consumption ($\eta_{jp}$) and production
         ($\mu_{jp}$) coefficients.}
\label{tab:stoichiometry}
\end{table}

\subsection{Graph modeling and solution}

\begin{figure}
    \centering
    \includegraphics[width=\linewidth]{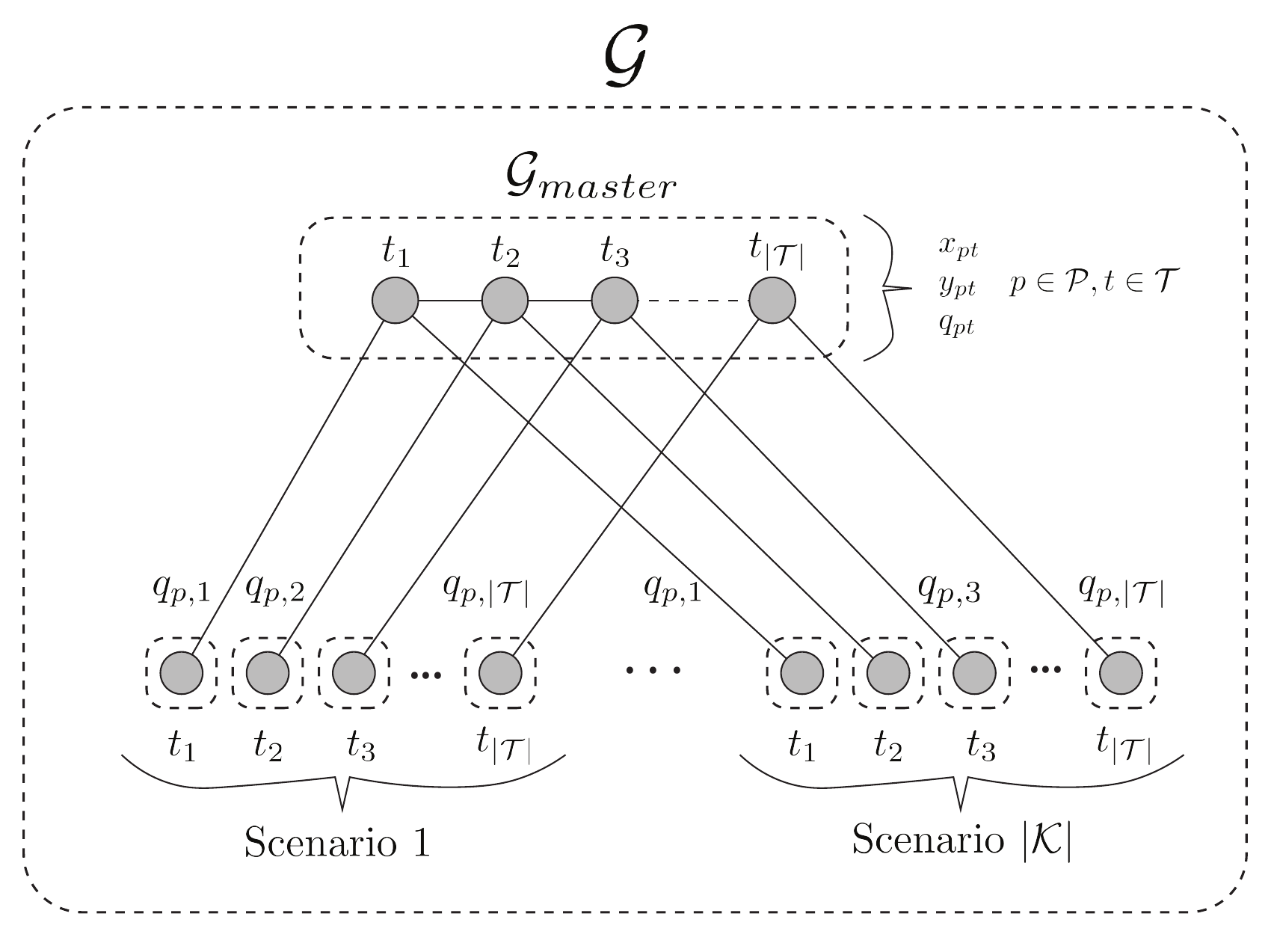}
    \caption{The graph structure of the capacity-expansion planning problem. The master problem subgraph, $\mathcal{G}_{master}$ contains the first-stage decisions across all time points. The second-stage capacity decisions, $q_{pt}$ are linked via linking constraints to the scenario subgraphs at each time point.}
    \label{fig:graph_case_study}
\end{figure}

The problem was implemented as a graph model within {\tt Plasmo.jl} and solved using BD within {\tt PlasmoBenders.jl}. The resulting OptiGraph is visualized in Figure \ref{fig:graph_case_study}. The first stage problem is placed on a subgraph called $\mathcal{G}_{master}$, which spans all time points (each time point represented as a node). The nodes of this subgraph contain the variables $x_{pt}$, $y_{pt}$, and $q_{pt}, p \in \mathcal{P}, t \in \mathcal{T}$ Each of these nodes is also linked to the previous time point (\eqref{eq:master_link}). The time points of each scenario are then added as independent subgraphs, each containing one node. These scenario time points are linked to the master problem by constraints (edges) containing the complicating variable $q_{pt}$. This results in $|\mathcal{K}| \times |\mathcal{T}|$ subproblem subgraphs that can be solved independently and in parallel. 

The graph model was solved using the mp-data and {\tt PlasmoBenders.jl}. To do this, the mp problem for this case study was solved once using the PPOPT toolbox \cite{kenefake2022ppopt}. The graph structure defined above was passed to {\tt PlasmoBenders.jl}, which automatically detects the subproblems once the master problem is set. A classic MILP solver (HiGHS or Gurobi) was set on the master problem subgraph, and the custom solver was set on each of the subproblem subgraphs. The resulting mp data was set within the custom solver used on each subproblem subgraph, and Benders was applied to solve each problem. 

To analyze the behavior of mp-Benders, we tested the framework using both HiGHS \cite{huangfu2018parallelizing} and Gurobi \cite{gurobi} as the MILP solver used for solving the master problem and the custom mp solver on the subproblems. We also solved the subproblems using HiGHS and Gurobi for comparison with mp. All tests were run using {\tt PlasmoBenders.jl} for consistency. Solves were done with both multi-cut and single-cut Benders approaches.

\section{Results and Discussion}
\label{sec:Results and Discussion}

\subsection{Computational Results: Multi-cut vs Single-cut Benders Decomposition}

\begin{figure}[!htp]
  \centering
  \begin{subfigure}[b]{0.45\textwidth}
    \centering
    \includegraphics[width=\textwidth]{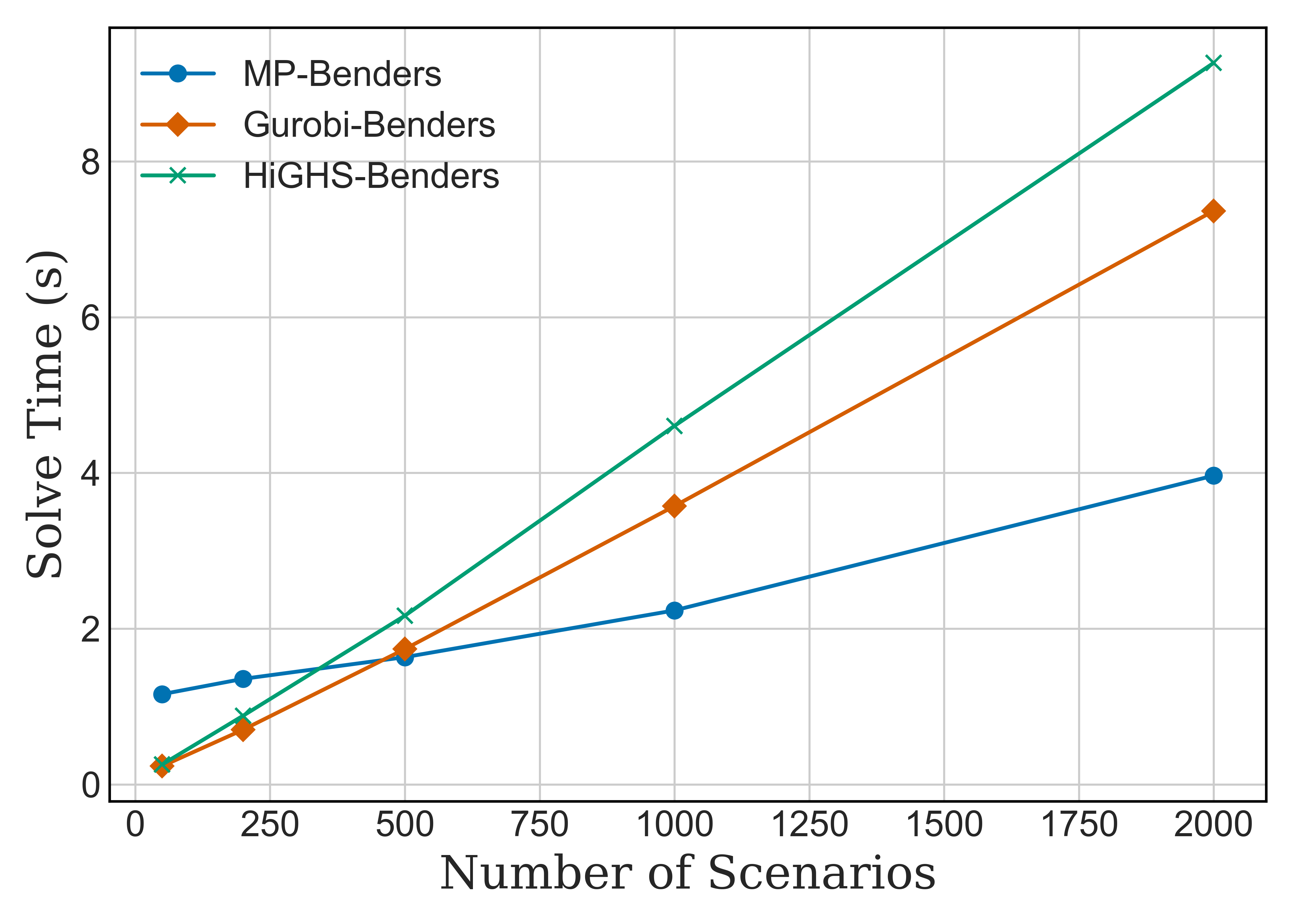}
    \caption{$T=5$}
    \label{fig:mc_5}
  \end{subfigure}
  \hfill
  \begin{subfigure}[b]{0.45\textwidth}
    \centering
    \includegraphics[width=\textwidth]{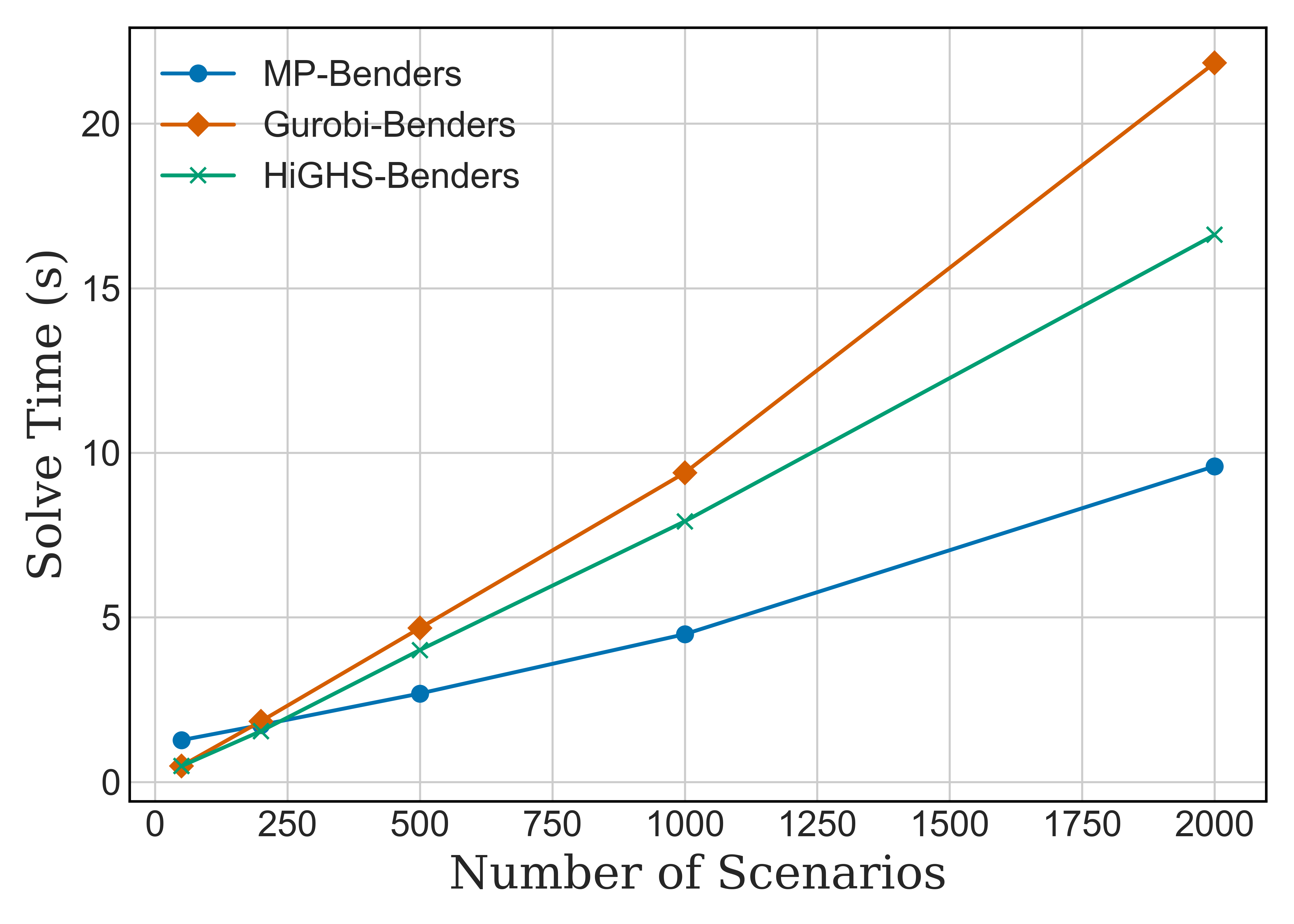}
    \caption{$T=10$}
    \label{fig:mc_10}
  \end{subfigure}

  \vspace{0.5cm}

  \begin{subfigure}[b]{0.45\textwidth}
    \centering
    \includegraphics[width=\textwidth]{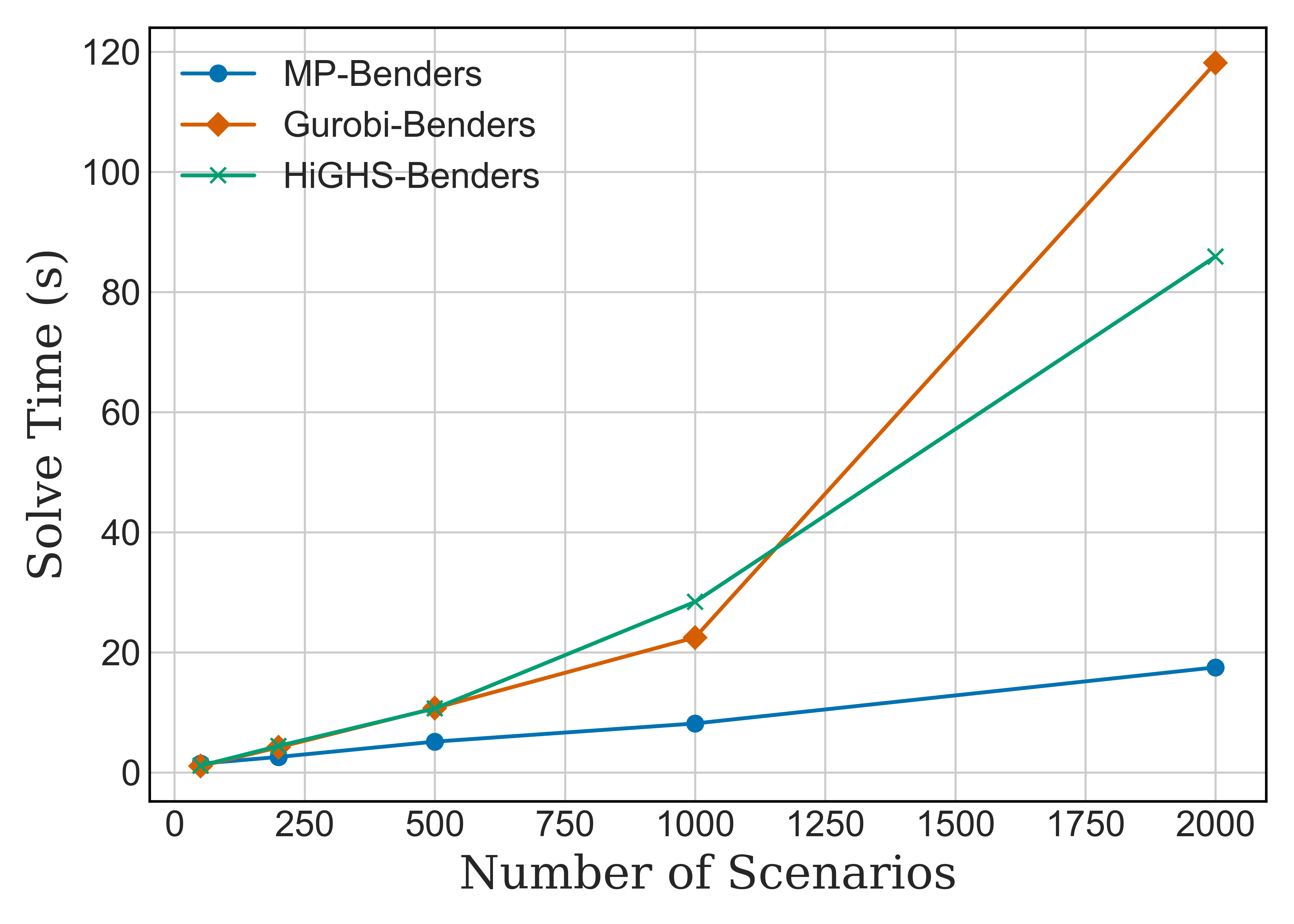}
    \caption{$T=25$}
    \label{fig:mc_25}
  \end{subfigure}
  \hfill
  \begin{subfigure}[b]{0.45\textwidth}
    \centering
    \includegraphics[width=\textwidth]{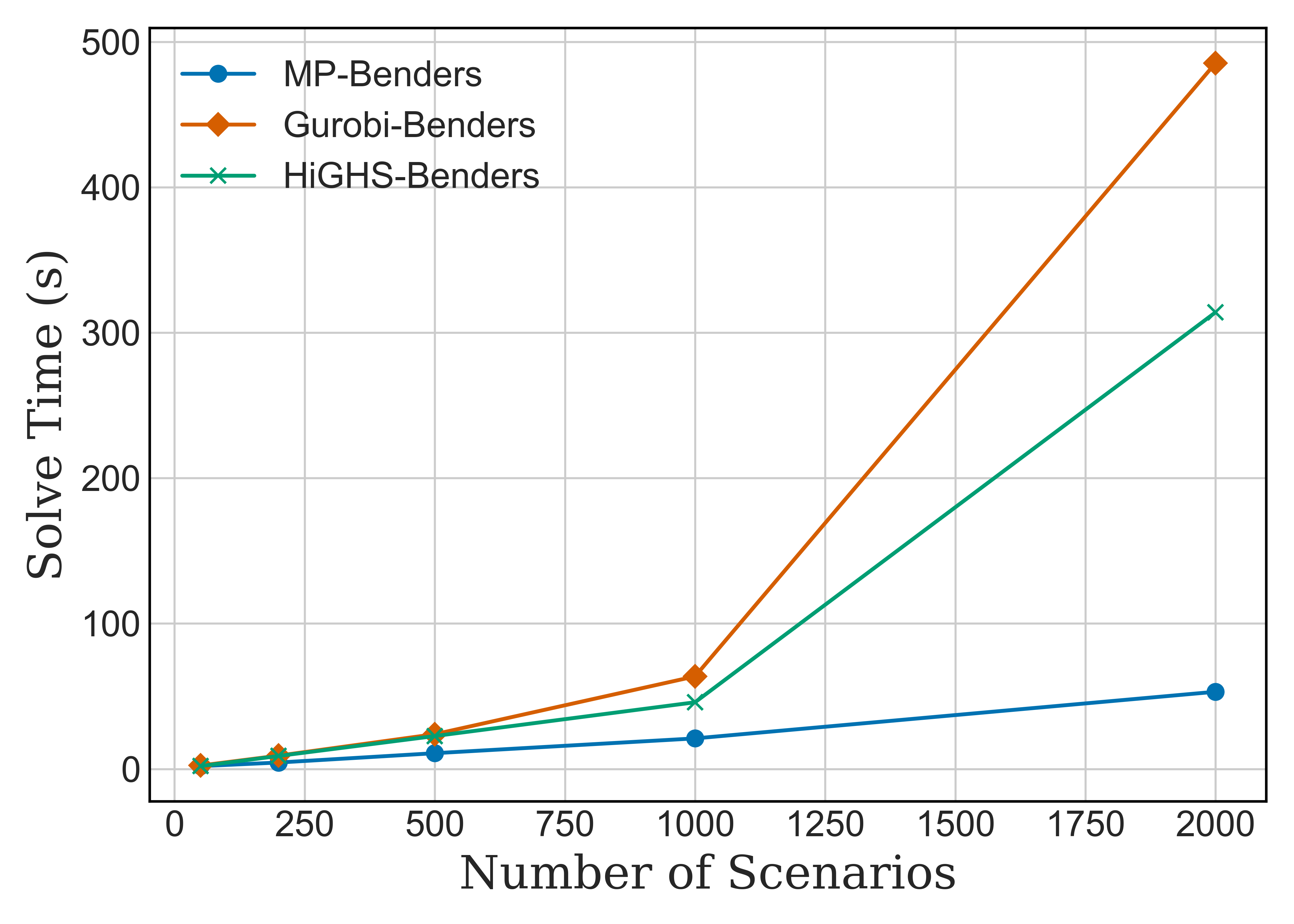}
    \caption{$T=50$}
    \label{fig:mc_50}
  \end{subfigure}

  \caption{Subproblem solve time vs. number of scenarios for different time horizons $T$ (multi-cut).}
  \label{fig:mc_all}
\end{figure}

Figures~\ref{fig:mc_5}--\ref{fig:mc_50} illustrate the computational results for the multi-cut Benders decomposition with planning horizons $T = 5, 10, 25, 50$. The subproblem solution time is reported for varying numbers of scenarios $K = 50, 200, 500, 1000, 2000$, comparing the mp-based critical region querying method (``mp-Benders'') to direct solves using Gurobi (``Gurobi-Benders'') and HiGHS (``HiGHS-Benders''). Across all cases, the mp-Benders approach consistently yields faster subproblem solution times, particularly as the number of scenarios increases. These performance improvements are primarily due to the avoidance of repeated LP solves by querying precomputed affine solutions from critical regions. It is important to note that the reported subproblem solution times for the mp-Benders case include the one-time overhead of computing the critical regions (offline). As a result, for smaller scenario sizes, the mp-Benders approach may not outperform the direct Gurobi or HiGHS subproblem solves, since the relative impact of the critical region generation time is more significant in such cases.

\begin{figure}[!htp]
  \centering
  \begin{subfigure}[b]{0.45\textwidth}
    \centering
    \includegraphics[width=\textwidth]{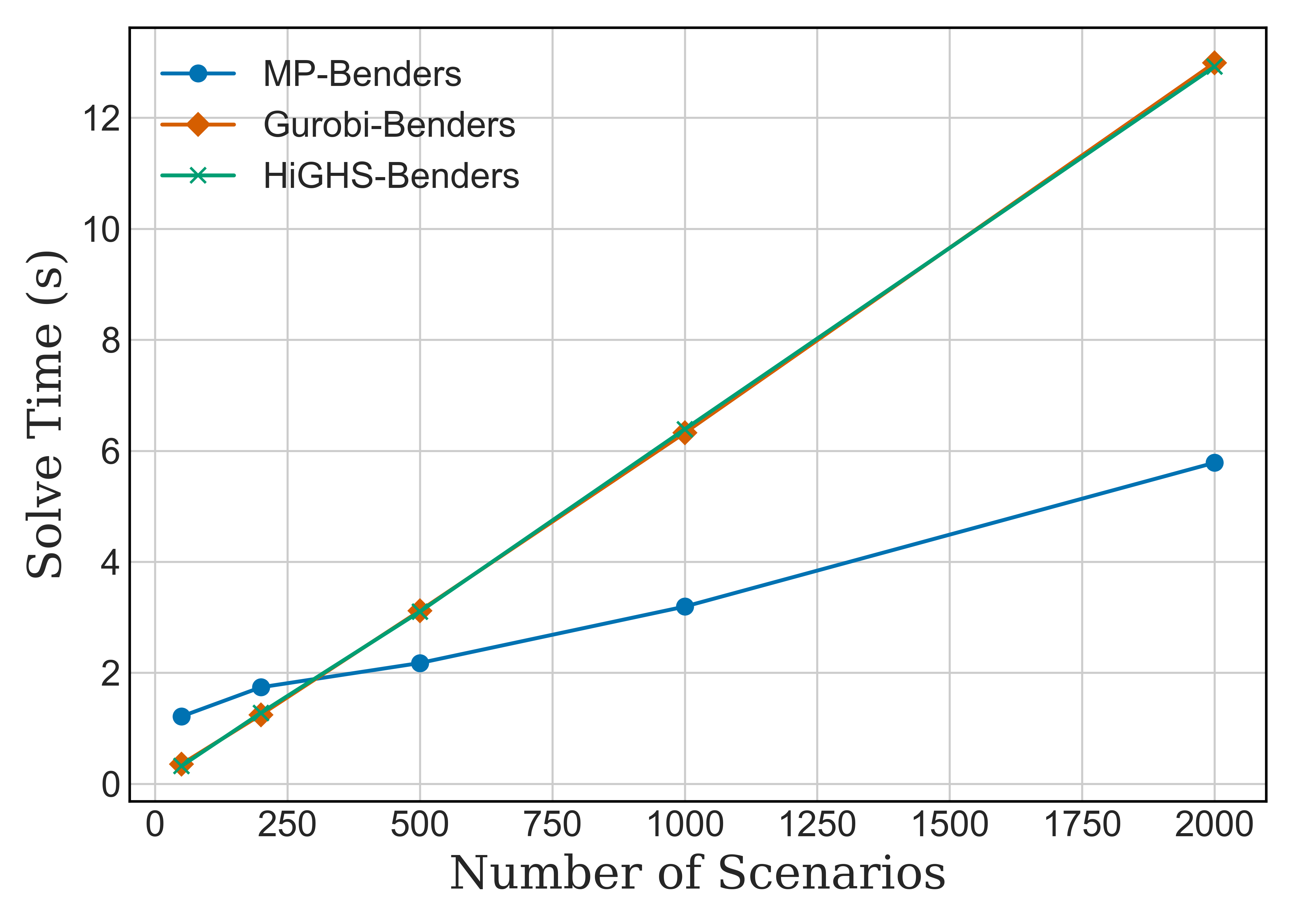}
    \caption{$T=5$}
    \label{fig:sc_5}
  \end{subfigure}
  \hfill
  \begin{subfigure}[b]{0.45\textwidth}
    \centering
    \includegraphics[width=\textwidth]{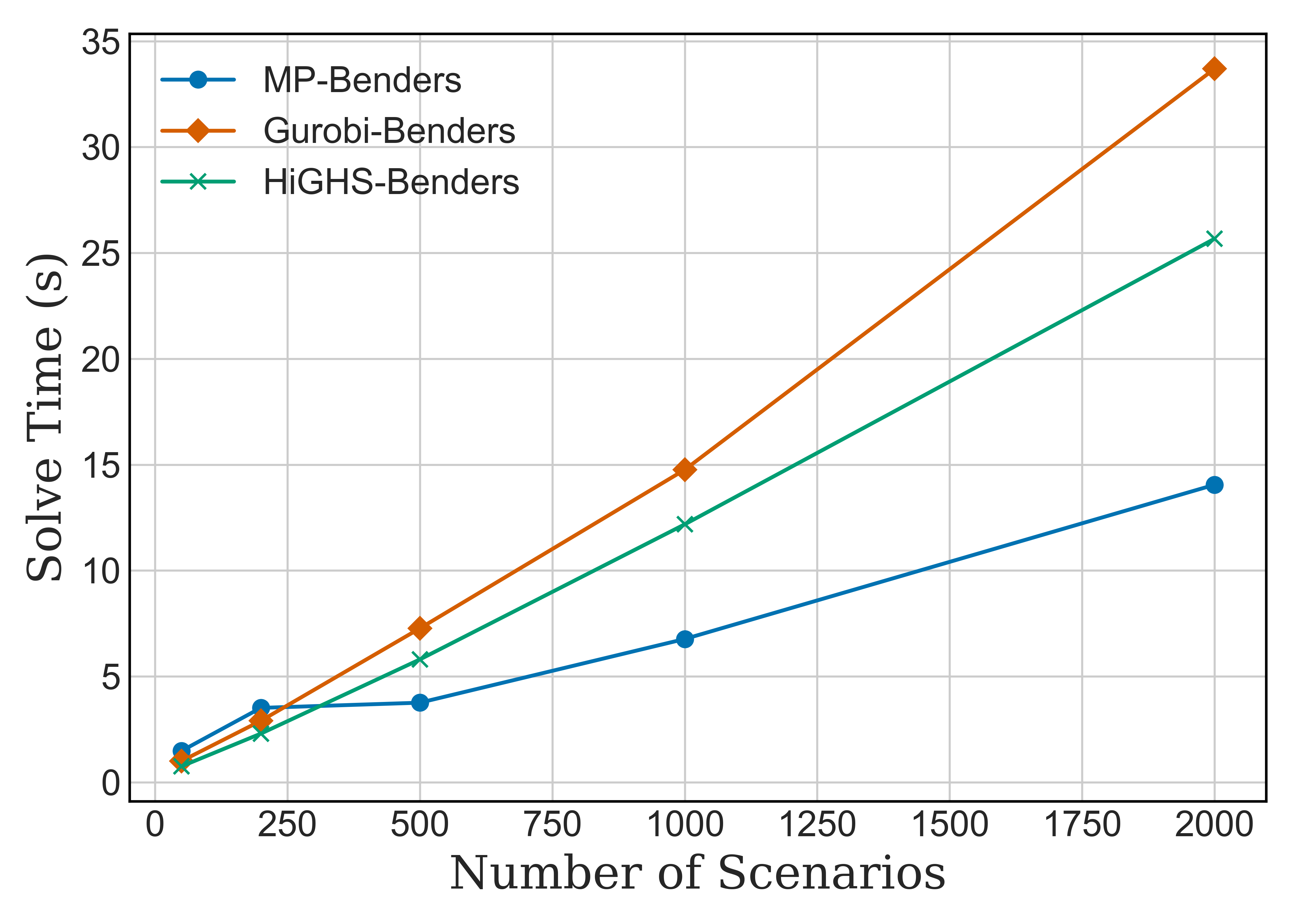}
    \caption{$T=10$}
    \label{fig:sc_10}
  \end{subfigure}

  \vspace{0.5cm}

  \begin{subfigure}[b]{0.45\textwidth}
    \centering
    \includegraphics[width=\textwidth]{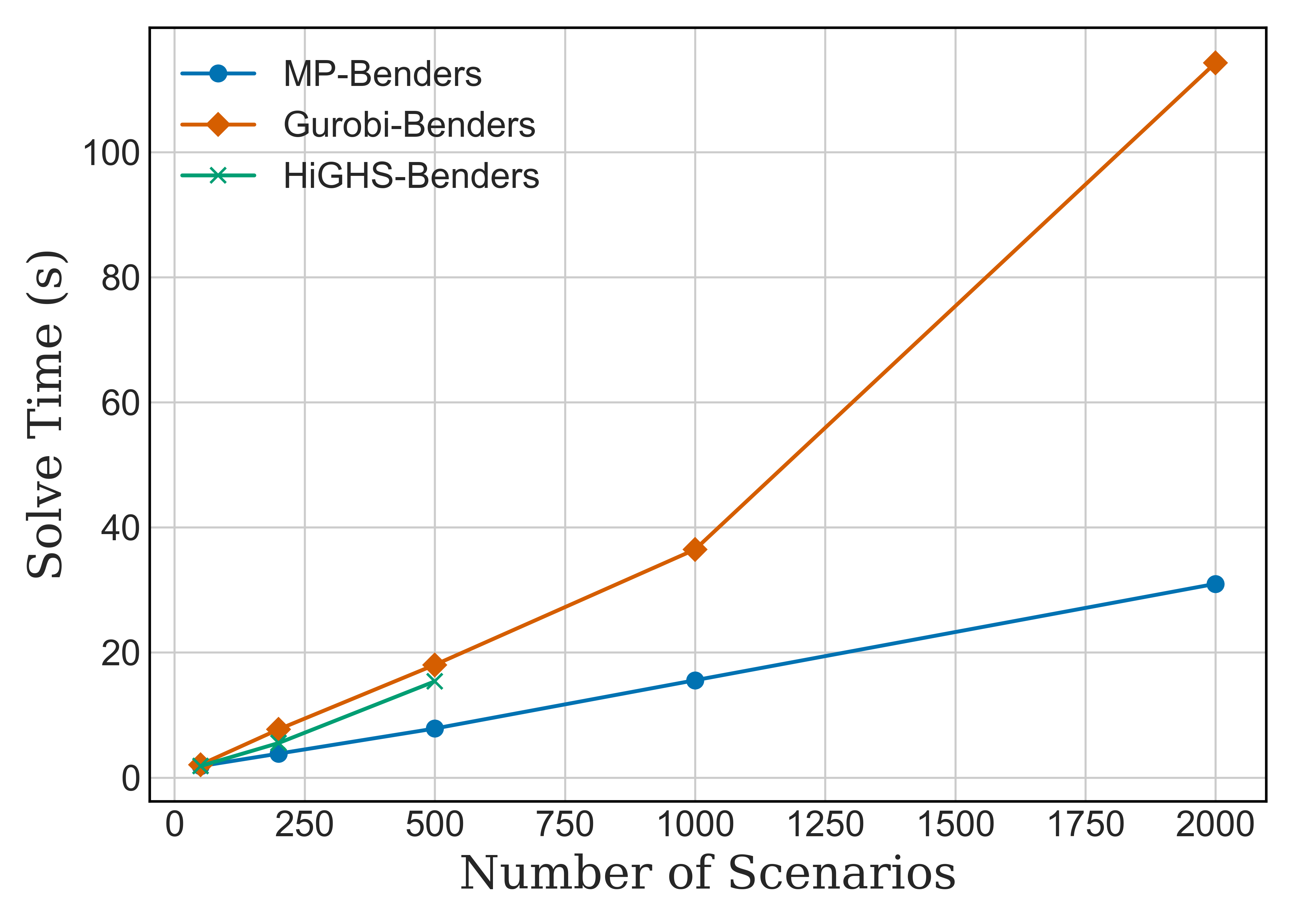}
    \caption{$T=25$}
    \label{fig:sc_25}
  \end{subfigure}
  \hfill
  \begin{subfigure}[b]{0.45\textwidth}
    \centering
    \includegraphics[width=\textwidth]{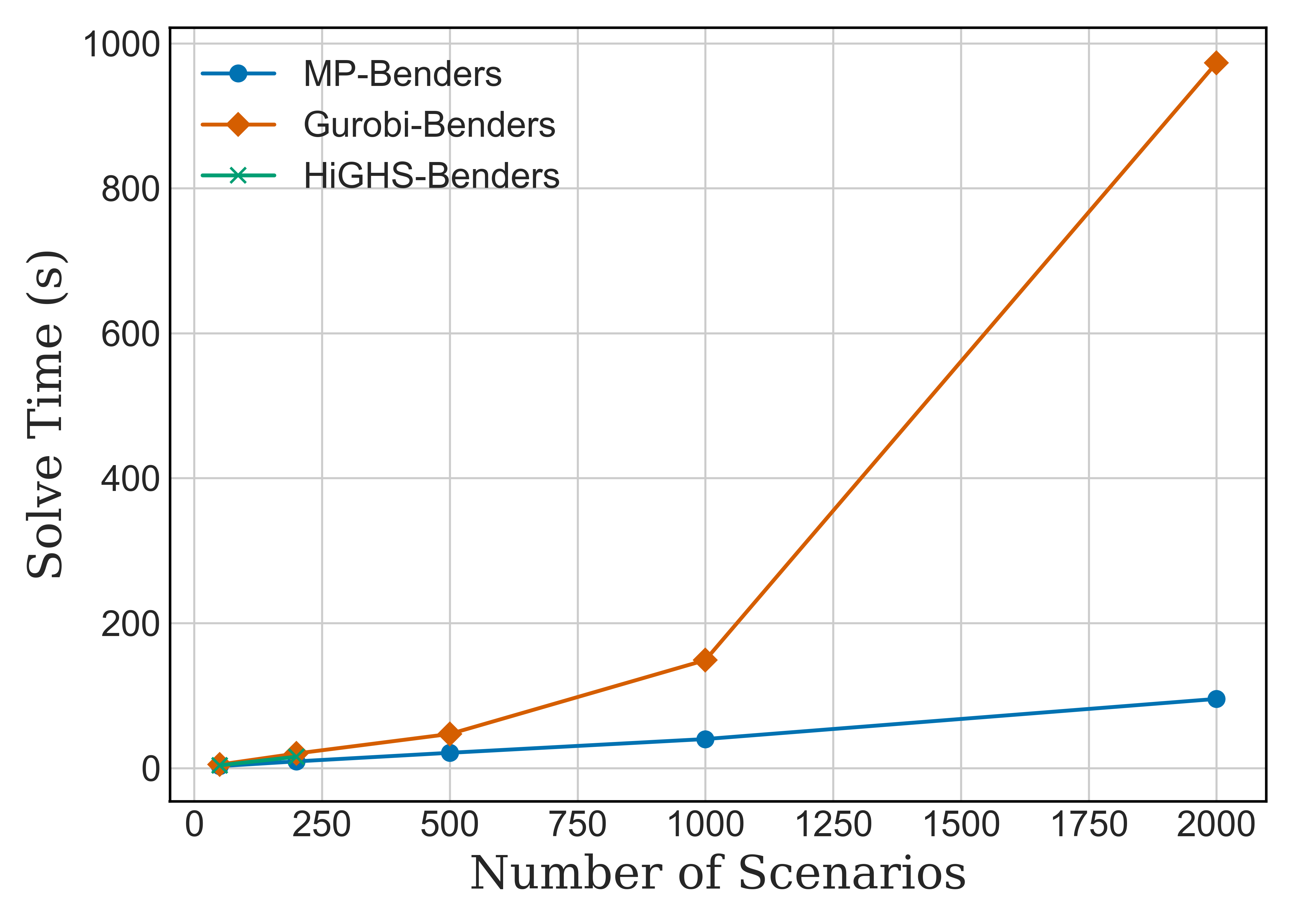}
    \caption{$T=50$}
    \label{fig:sc_50}
  \end{subfigure}

  \caption{Subproblem solve time vs. number of scenarios for different time horizons $T$ (single-cut).}
  \label{fig:sc_all}
\end{figure}

Figures~\ref{fig:sc_5}--\ref{fig:sc_50} present the results for the single-cut Benders decomposition variant. The primary motivation for examining this variant is to assess the performance advantages of mp-Benders surrogates in settings that require a large number of Benders iterations, as is typical with single-cut formulations. Notably, the subproblem solve time using mp-Benders is significantly lower than that of Gurobi- or HiGHS solvers, and this difference is pronounced compared to the multi-cut case. This observation suggests that the computational benefits of mp surrogates become increasingly prominent in scenarios where the decomposition algorithm must generate a greater number of cuts to converge.

Nevertheless, for this particular case study, all single-cut variants—mp-Benders, Gurobi-Benders, and HiGHS-Benders exhibit longer overall solution times relative to their multi-cut counterparts. This outcome is primarily attributed to the slower convergence of the single-cut scheme, which aggregates recourse information and thus provides less informative feedback per iteration.

It is also important to highlight that for the large planning horizons ($T = 25$ and $T = 50$), the HiGHS-Benders implementation failed to solve problem instances with high scenario counts ($K \in \{500, 1000, 2000\}$). Despite multiple attempts, we were unable to isolate the root cause of these failures, which we suspect stem from numerical instability or solver tolerance thresholds. Consequently, the affected data points are omitted from the corresponding plots \ref{fig:sc_25} and \ref{fig:sc_50}.

\subsection{Interpretability of Benders Decomposition using mp programming}

Beyond computational gains, one of the most compelling benefits of the mp-based surrogate framework is its ability to enhance the analysis and interpretability of algorithmic and solution behavior. Specifically, stochastic programming models are powerful, but their solutions are often difficult for users to understand and therefore trust, which has hindered their adoption in practice \cite{rathi2024enhancing}. A key challenge is moving beyond simply obtaining an optimal solution to explaining \textit{why} that particular solution is optimal \cite{bertsimas2021voice}. Classical Benders decomposition often exacerbates this issue by treating the subproblem solver as a black box, offering little transparency into how uncertainty and first-stage decisions interact with the recourse solutions.

To address this challenge, researchers such as  Rathi \textit{et al.} \cite{rathi2024enhancing} have proposed systematic methods to enhance the explainability of SP solutions. Their approach focuses on creating reduced, more interpretable models that capture the key features of the original problem. They achieve this through two post-hoc techniques: (i) \textbf{Scenario Reduction}, which clusters the vast number of scenarios based on the similarity of their optimal \textit{recourse decisions} \cite{bounitsis2022data} rather than just their input data, thereby identifying the features of uncertainty that are most relevant to the solution; and (ii) \textbf{Recourse Reduction}, which identifies a small set of ``principal recourse variables'' that are most critical for achieving an optimal or near-optimal solution, simplifying the problem structure. The explicit goal is to create a simpler model that allows a user to better understand why the optimal first-stage decisions were made.

The proposed mp framework provides a direct and analytical way to achieve interpretability. The mp formulation provides an explicit, geometric partitioning of the parameter space into critical regions, where each region is associated with a unique affine decision and value function. These critical regions can be viewed as an analytical form of the recourse-based scenario clustering proposed by Rathi \textit{et al.} \cite{rathi2024enhancing}; all scenarios (parameter vectors) that fall within a single critical region are, by definition, a ``cluster'' that shares the same set of active constraints and is governed by the same affine recourse function. This allows us to understand entire classes of scenarios at once. Furthermore, by inspecting the affine solution function ($x^{*} = A^{v}\theta + b^{v}$) within each region, we can perform an analysis analogous to recourse reduction. The components of the matrix $A^v$ reveal which recourse variables are sensitive to changes in first-stage decisions and uncertain parameters, thereby identifying the most ``principal'' or flexible decisions for that entire region.

By embedding this structure into the Benders framework, we gain the ability to visualize how the master problem iterates navigate this explicit critical region landscape. We have shown this for one time-scenario combination $t=5$ and $K=500$. Figures 7-10 depict decision trajectories through these regions for the multi-cut and single-cut
schemes, respectively, across time steps t = 1 to t = 5. This visualization reveals several key insights:

\begin{itemize}
    \item \textbf{Convergence behavior:} In the multi-cut case, the iterates tend to converge more directly toward the final solution, traversing fewer critical regions. This aligns with the theoretical observation that multi-cut Benders converges in fewer iterations due to richer feedback from the subproblems.

    \item \textbf{Decision consistency:} Within certain time steps (e.g., $t=3$ or $t=4$), the iterates remain confined to a single or a small number of regions. This suggests robustness of the recourse structure with respect to small perturbations in the master decision, i.e., the recourse decisions exhibit local linear stability.

    \item \textbf{Aggregation effect:} The single-cut formulation requires more iterations and exhibits more oscillatory behavior, with the decision path traversing a greater number of distinct regions. This reflects the loss of resolution caused by aggregating all subproblem feedback into a single surrogate.

    \item \textbf{Policy transparency:} Because each region corresponds to a fixed active set of constraints, we can determine \textit{why} a particular recourse policy (the affine function) is optimal for an entire class of scenarios (the critical region) by identifying the binding operational constraints. This allows for a richer post-optimal analysis, identifying which constraints dominate solution behavior under different uncertainty realizations.
\end{itemize}

Taken together, these insights offer a novel perspective on stochastic optimization algorithms: not just as black-box iterative solvers, but as structured geometric processes navigating a known (albeit high-dimensional) partition of the parameter space. This geometric lens can support better debugging, policy interpretation, and even sensitivity analysis with respect to scenario inputs or problem data, directly addressing the need for greater explainability in stochastic optimization.

\newpage

\begin{figure}[!htp]
  \centering
  \begin{subfigure}[b]{0.45\textwidth}
    \centering
    \includegraphics[width=\textwidth]{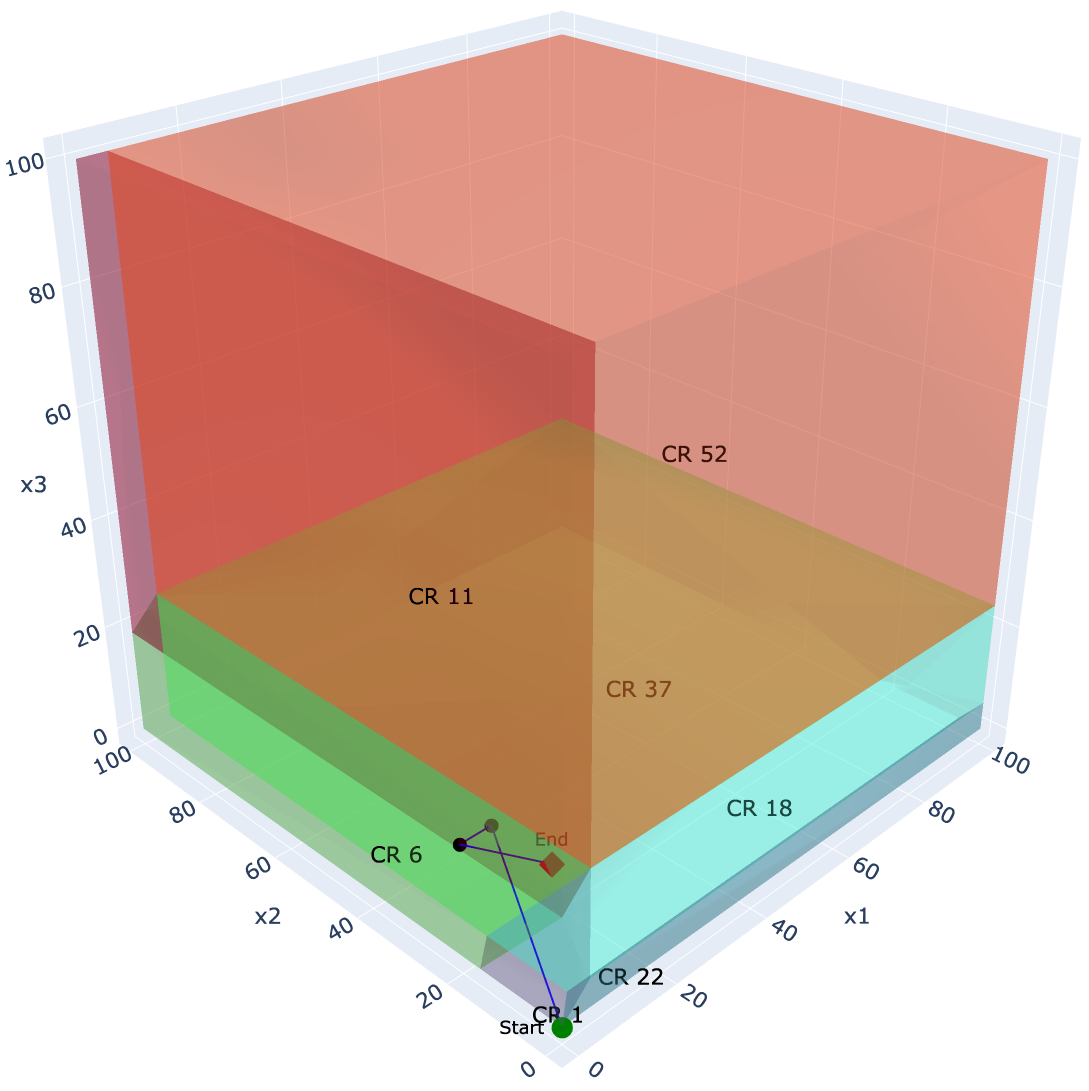}
    \caption{$t=1$}
    \label{fig:mc_polytope_1}
  \end{subfigure}
  \begin{subfigure}[b]{0.45\textwidth}
    \centering
    \includegraphics[width=\textwidth]{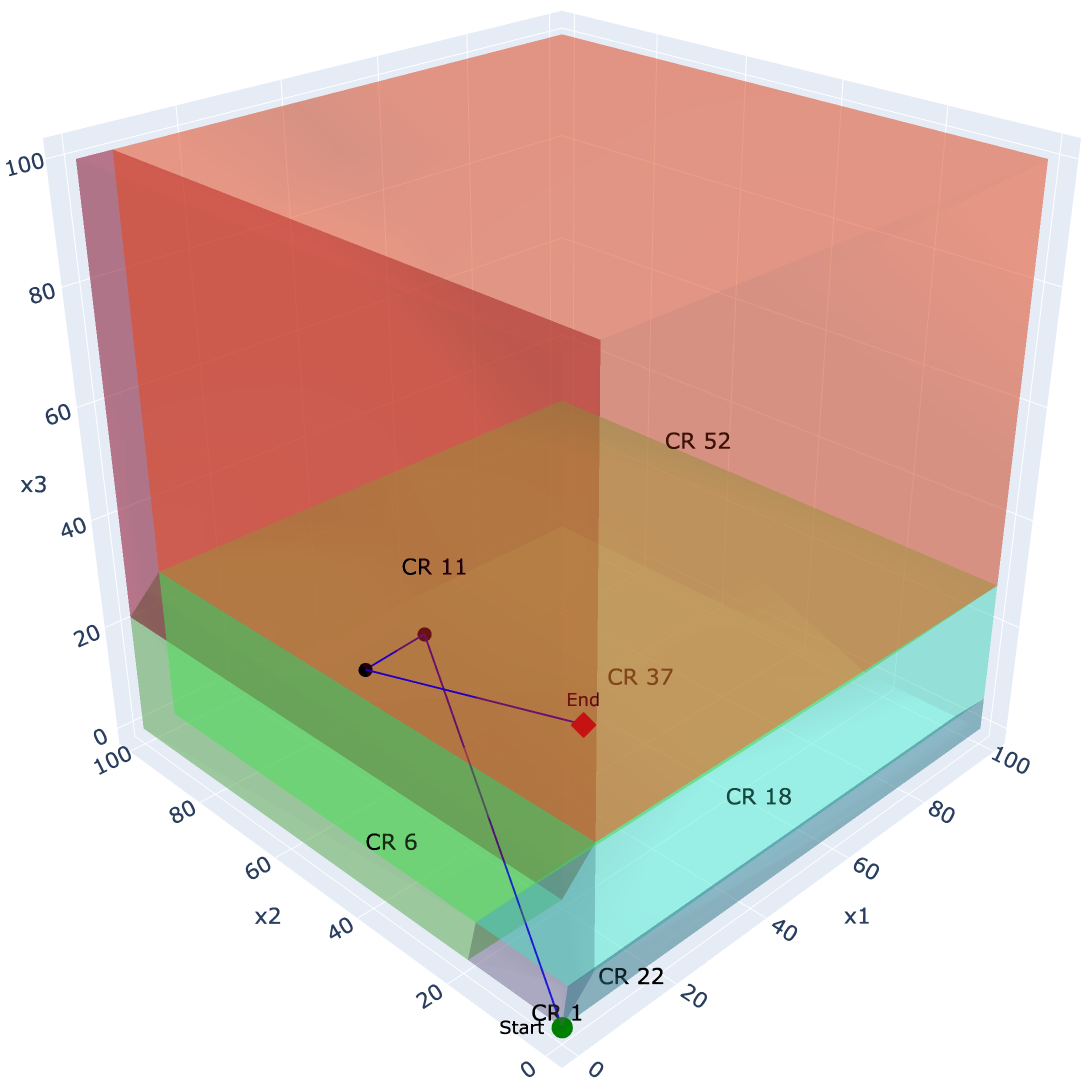}
    \caption{$t=2$}
    \label{fig:mc_polytope_2}
  \end{subfigure}
  \vspace{0.3cm}
  \begin{subfigure}[b]{0.45\textwidth}
    \centering
    \includegraphics[width=\textwidth]{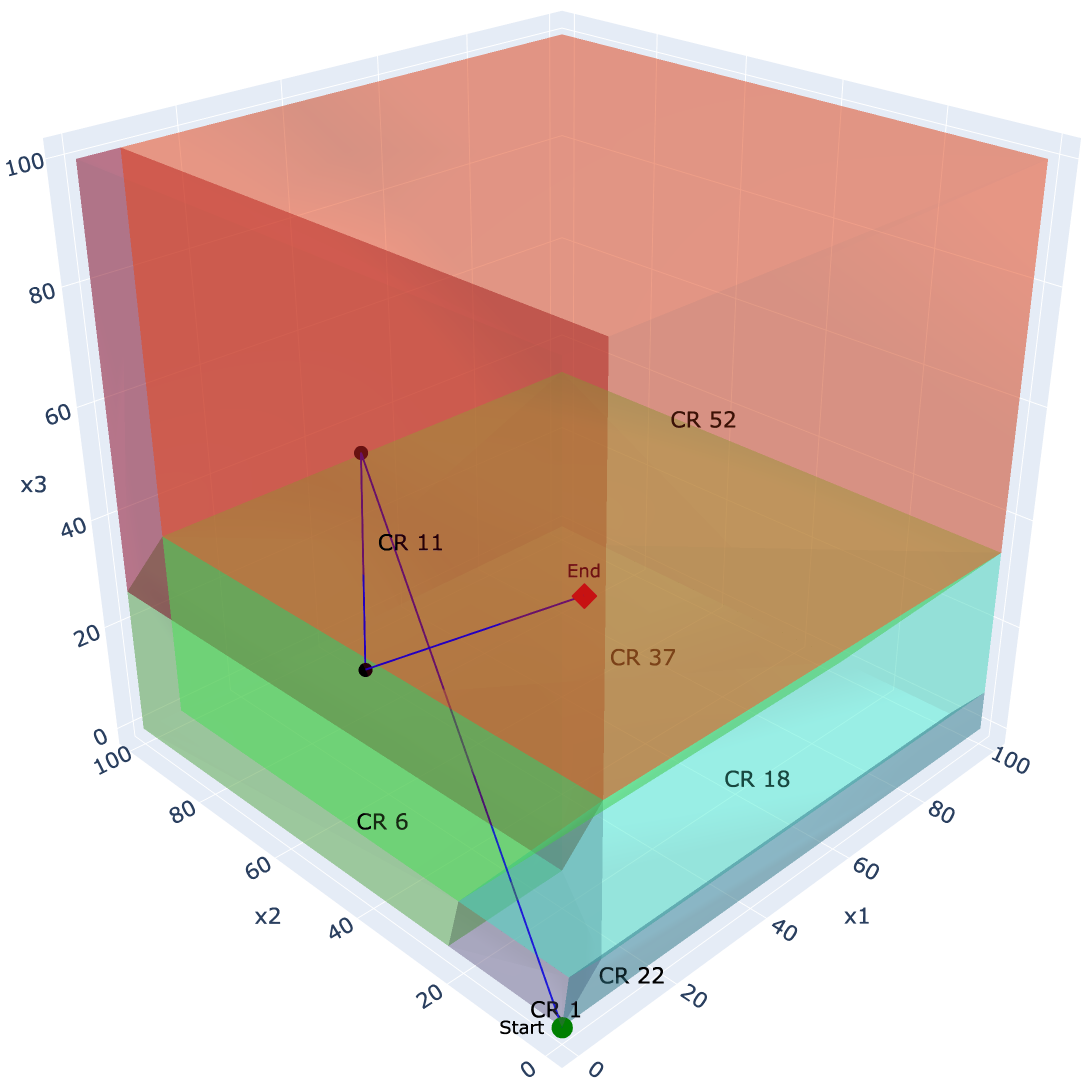}
    \caption{$t=3$}
    \label{fig:mc_polytope_3}
  \end{subfigure}
  \begin{subfigure}[b]{0.45\textwidth}
    \centering
    \includegraphics[width=\textwidth]{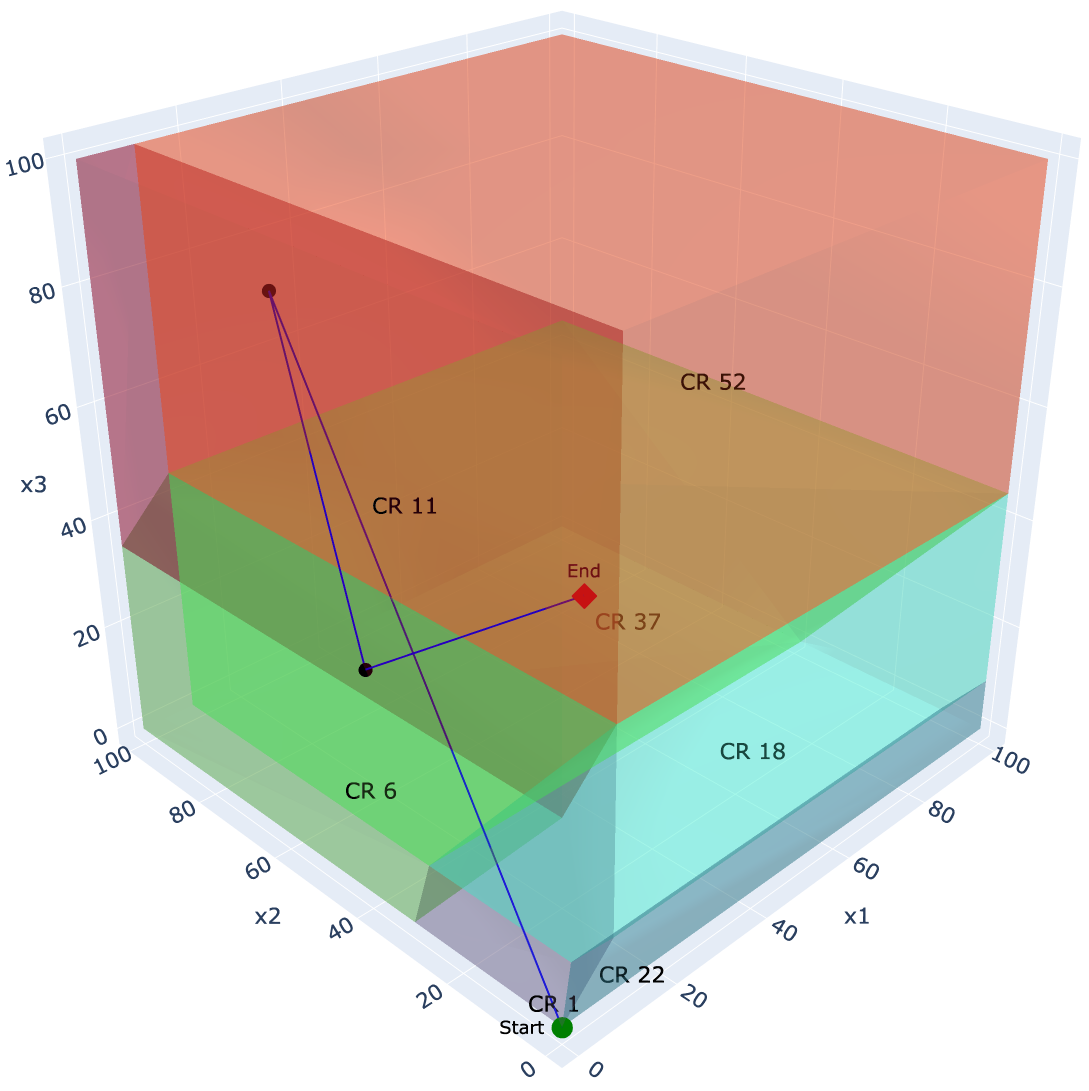}
    \caption{$t=4$}
    \label{fig:mc_polytope_4}
  \end{subfigure}

  \caption{Multi-cut: Decision trajectory in critical region space for $t=1$ to $t=4$.}
  \label{fig:mc_polytope_part1}
\end{figure}

\newpage

\begin{figure}[!htp]
  \centering
  \begin{subfigure}[b]{0.9\textwidth}
    \centering
    \includegraphics[width=\textwidth]{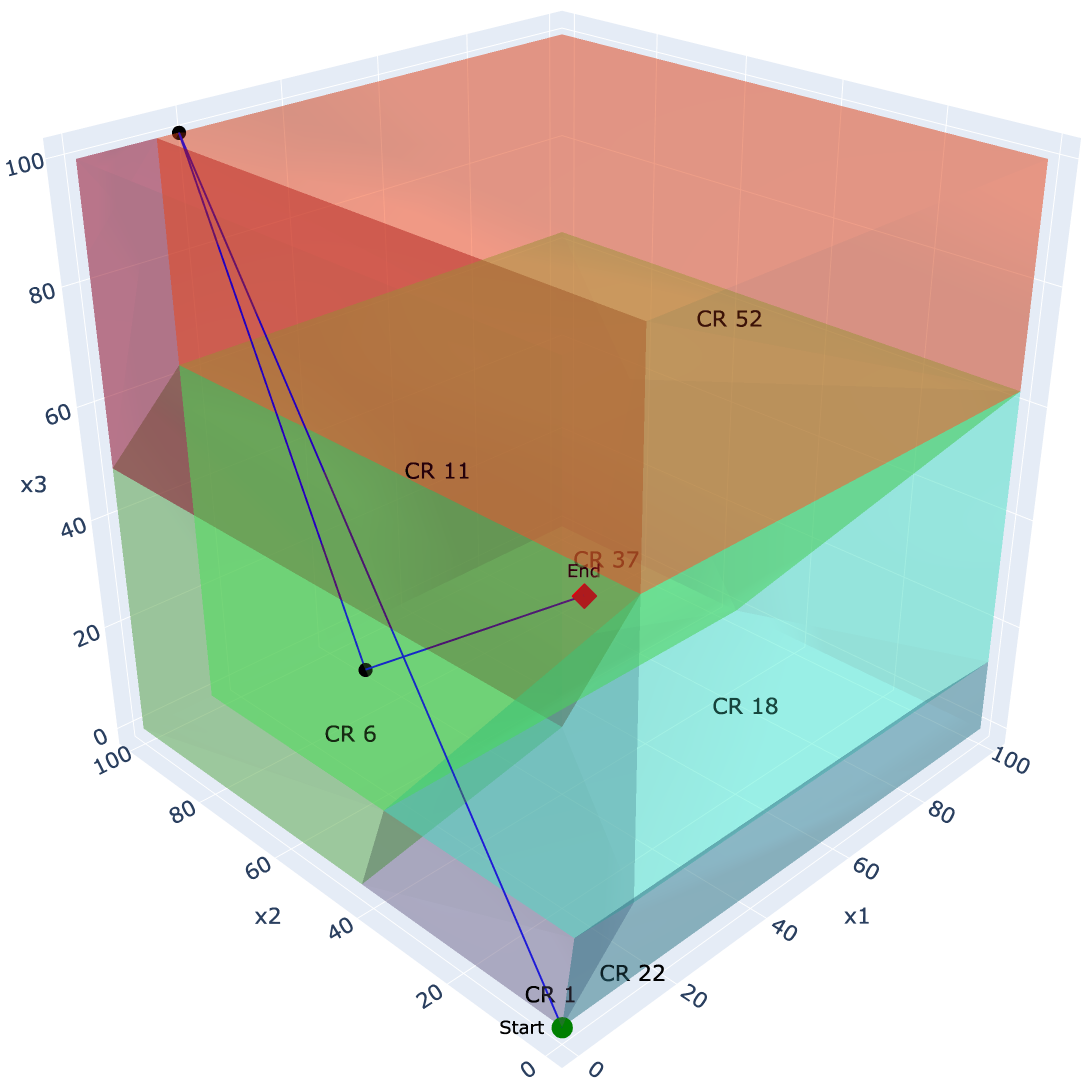}
    \caption{$t=5$}
    \label{fig:mc_polytope_5}
  \end{subfigure}

  \caption{Multi-cut: Decision trajectory in critical region space for $t=5$.}
  \label{fig:mc_polytope_part2}
\end{figure}

\newpage


\begin{figure}[!htp]
  \centering
  \begin{subfigure}[b]{0.45\textwidth}
    \centering
    \includegraphics[width=\textwidth]{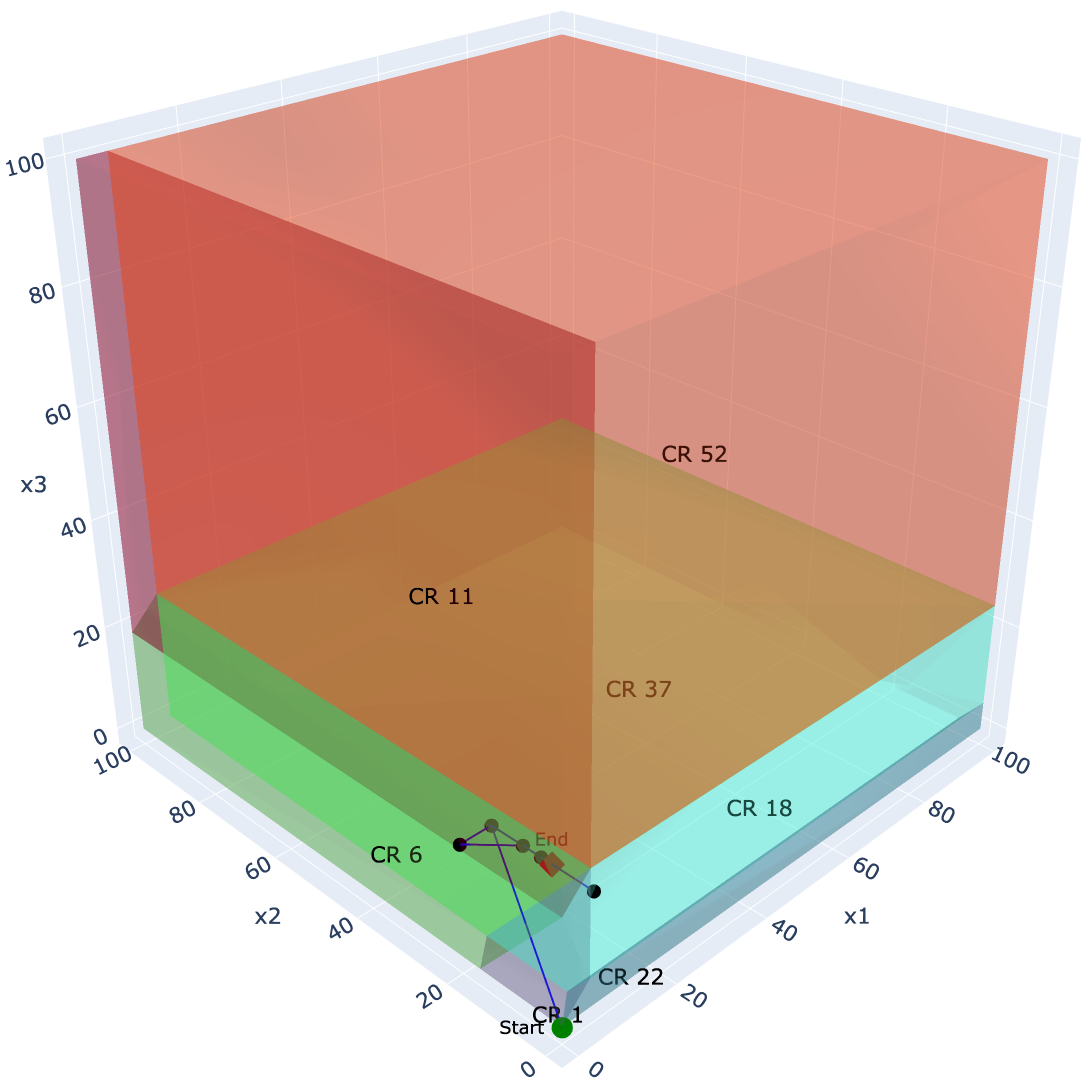}
    \caption{$t=1$}
    \label{fig:sc_polytope_1}
  \end{subfigure}
  \begin{subfigure}[b]{0.45\textwidth}
    \centering
    \includegraphics[width=\textwidth]{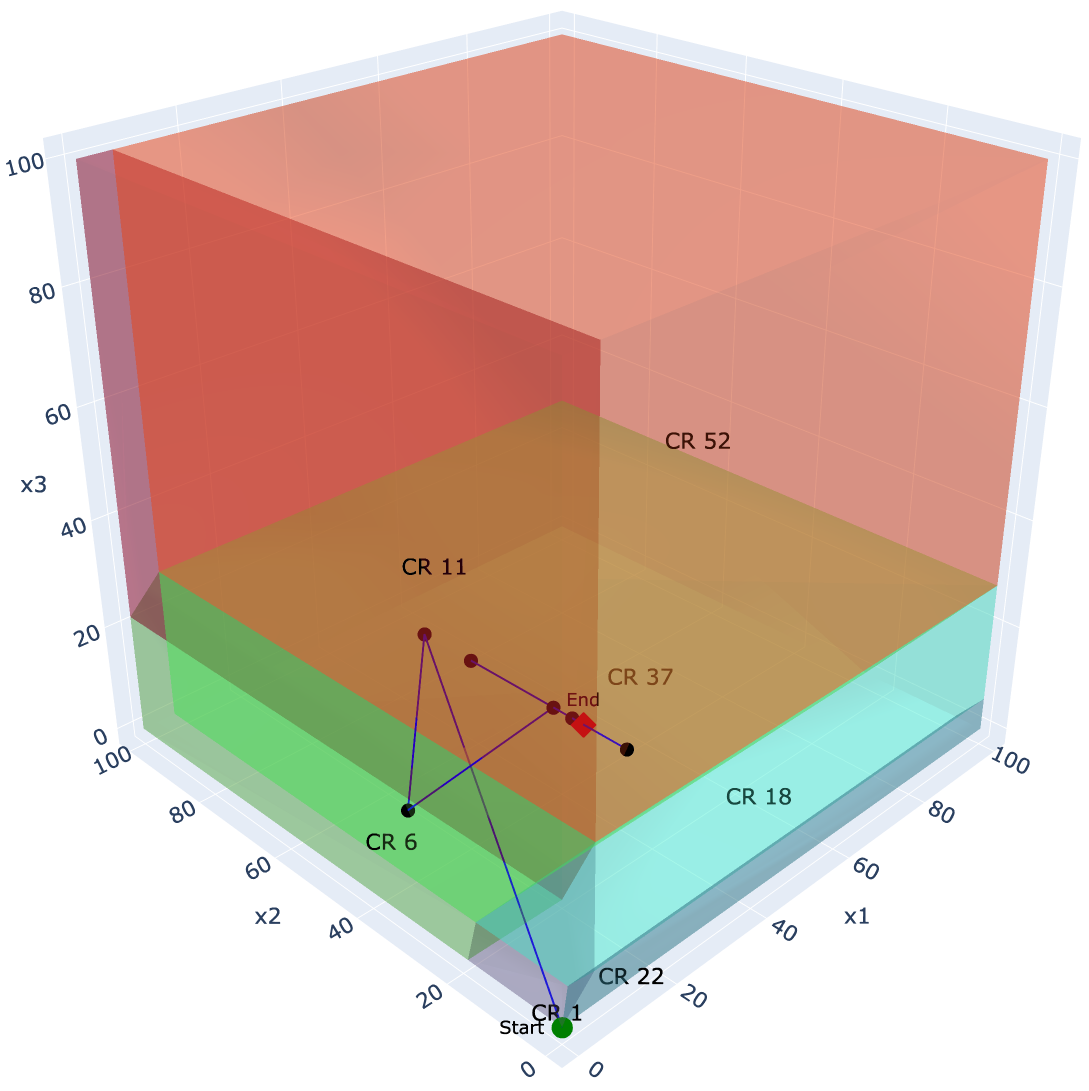}
    \caption{$t=2$}
    \label{fig:sc_polytope_2}
  \end{subfigure}

  \vspace{0.3cm}

  \begin{subfigure}[b]{0.45\textwidth}
    \centering
    \includegraphics[width=\textwidth]{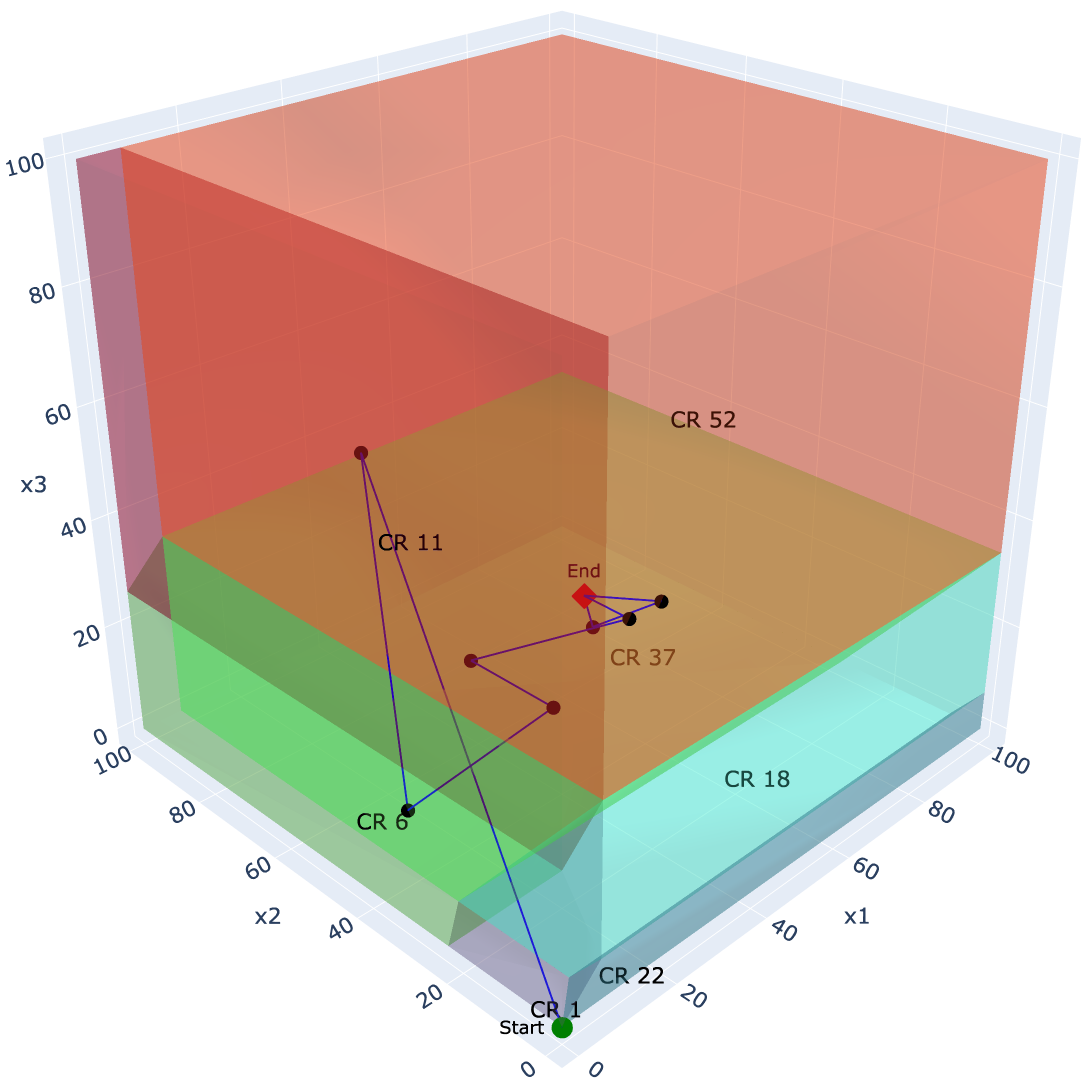}
    \caption{$t=3$}
    \label{fig:sc_polytope_3}
  \end{subfigure}
  \begin{subfigure}[b]{0.45\textwidth}
    \centering
    \includegraphics[width=\textwidth]{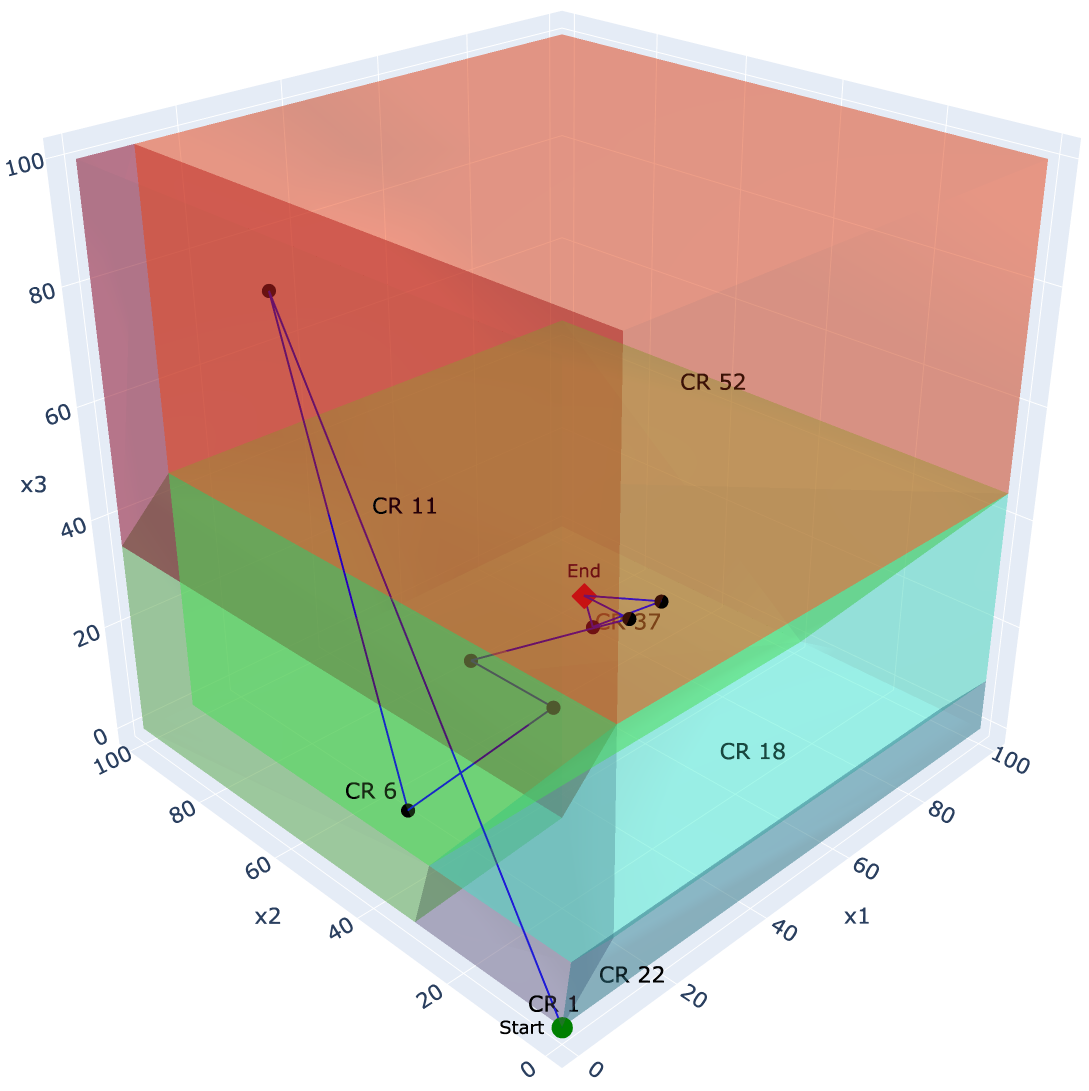}
    \caption{$t=4$}
    \label{fig:sc_polytope_4}
  \end{subfigure}

  \caption{Single-cut: Decision trajectory in critical region space for time steps $t=1$ to $t=4$.}
  \label{fig:sc_polytope_part1}
\end{figure}

\newpage

\begin{figure}[!htp]
  \centering
  \begin{subfigure}[b]{0.9\textwidth}
    \centering
    \includegraphics[width=\textwidth]{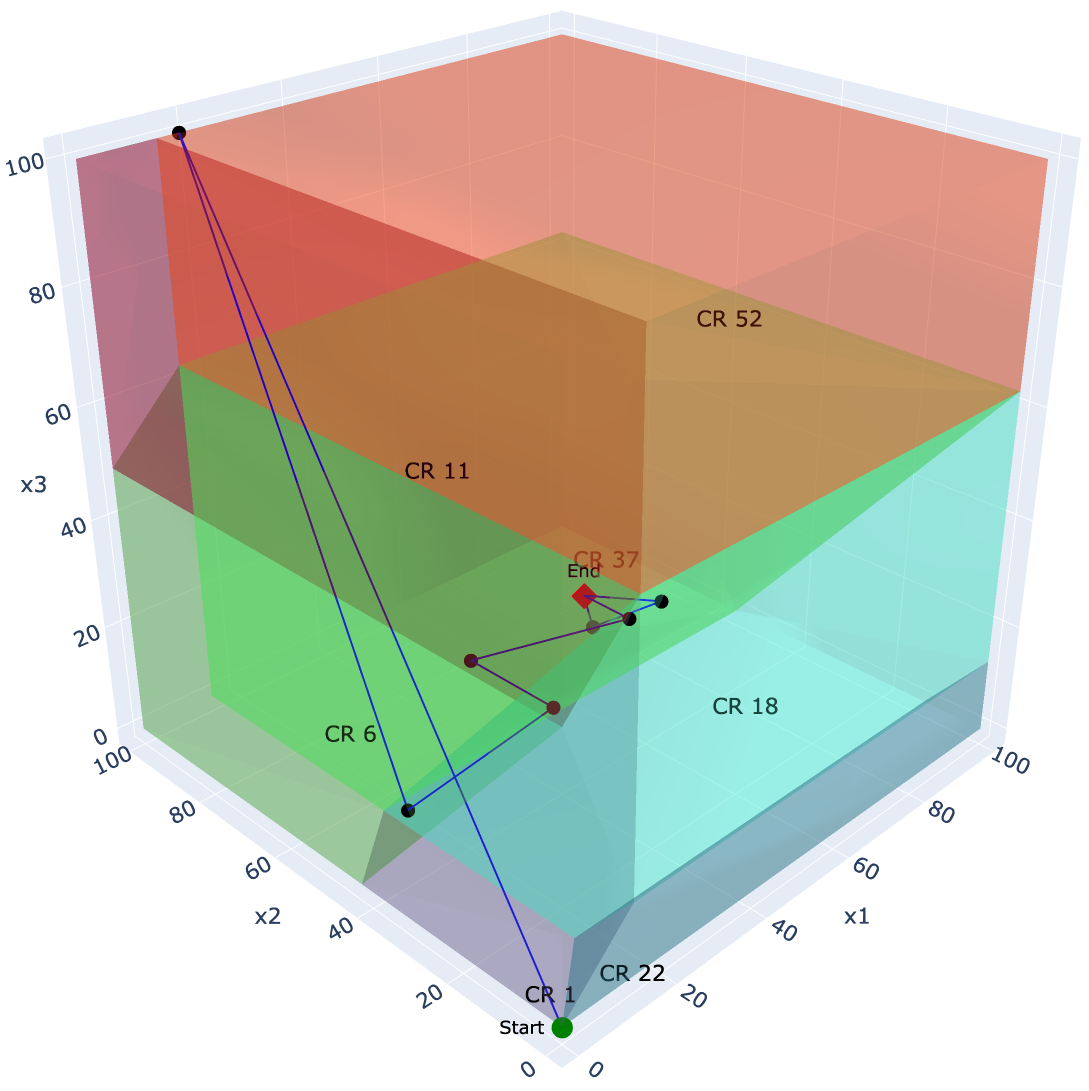}
    \caption{$t=5$}
    \label{fig:sc_polytope_5}
  \end{subfigure}

  \caption{Single-cut: Decision trajectory in critical region space for time step $t=5$.}
  \label{fig:sc_polytope_part2}
\end{figure}

\newpage

\begin{figure}[!htp]
    \centering
    \includegraphics[width=0.7\linewidth]{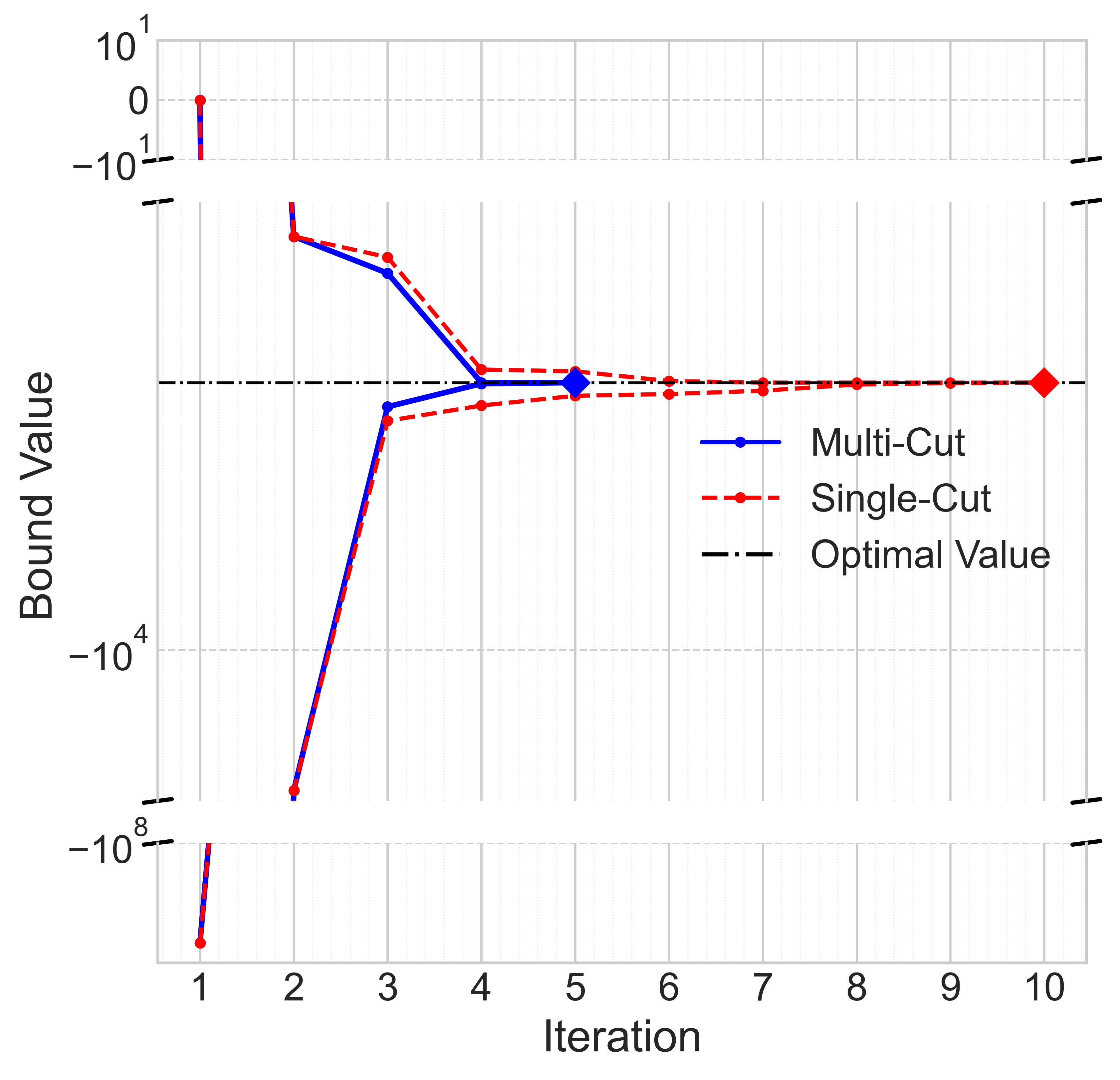}
    \caption{Bound convergence for multi-cut (blue) and single-cut (red) Benders algorithms with $T=5$ and $K=500$. Diamond markers indicate convergence to the optimal value (black dash-dot line). Note: The results shown are same for the mp-Benders method using Gurobi.}

    \label{fig:Upper_lower_bound}
\end{figure}

Figure \ref{fig:Upper_lower_bound} shows the evolution of the upper and lower bound for both multi-cut and single-cut cases for a representative problem instance.  This result confirms that the mp-based Benders approach converges to the optimal solution. This also illustrates that mp surrogates provide optimality guarantees, since primal and dual variables are optimal (not approximate).

\subsection{Benefits of Graph-Based Software Framework}

The graph modeling approach provides a flexible and modular framework for constructing and visualizing complex optimization problems. {\tt Plasmo.jl} allows for treating pieces of the optimization problem (subgraphs) as independent optimization problems, making it simpler to define decomposition algorithms that solve subproblems of a larger, overall problem. Furthermore, this allows for performing custom tasks on a single subgraph of the optimization problem without altering the overall problem structure or changing other subgraphs in the optimization problem. For instance, this could facilitate replacing a subgraph with a surrogate model or, as we did in this work, setting a custom solver on individual subgraphs that acted as an mp surrogate model.


\section{Conclusions and Future Work}
\label{sec:Conclusion}

We presented a computational framework for Benders decomposition that embeds multi-parametric programming (mp) surrogates within a graph representation of structured optimization problems. By integrating mp-LP surrogates of the graph node subproblems, we replace expensive repeated LP solves with region look-ups and affine evaluations, leading to substantial runtime reductions. Our implementation in {\tt PlasmoBenders.jl} demonstrates the flexibility of graph  decomposition and the interpretability benefits offered by critical region tracking. Computational experiments on a stochastic capacity expansion case study confirmed that the proposed mp-Benders approach provides consistent speedups across varying time horizons and scenario sizes, especially in the multi-cut variant.

We acknowledge a key limitation of the mp-based surrogate approach: its scalability with respect to subproblem size. As the number of uncertain parameters grows, the number of critical regions increases combinatorially, making enumeration and storage of mp solutions increasingly expensive. This poses challenges for problems with high-dimensional uncertainty or complex structure. As a direction for future work, we aim to develop more scalable, modular mp surrogates, where each subproblem may be approximated by multiple smaller mp models or decomposable sub-regions. This could mitigate the explosion in region count by partitioning the parameter space adaptively.

Furthermore, our current framework is centered around explicit parametric surrogates. We intend to extend the surrogate capabilities of the {\tt Plasmo.jl} and {\tt PlasmoBenders} ecosystem by incorporating alternative surrogate models from the machine learning domain such as decision trees, input convex neural networks (ICNNs), and other neural architectures as plug-and-play surrogates within the graph-based modeling approach. This will enable hybrid surrogate strategies that trade off interpretability, scalability, and approximation quality depending on the application context.

By building toward a unified, modular, and extensible surrogate-enhanced optimization framework, we hope to broaden the applicability of decomposition methods to real-world large-scale optimization problems.

\section*{Acknowledgments} 

This material is based on work supported by the U.S. Department of Energy under grant DE-0002722. We also acknowledge support from the National Science Foundation under award CMMI-2328160.

\section*{Nomenclature}

\begin{longtable}{ll}
\hline
\textbf{Symbol} & \textbf{Description} \\
\hline
\endfirsthead
\multicolumn{2}{c}%
{{\bfseries \tablename\ \thetable{} -- continued from previous page}} \\
\hline
\textbf{Symbol} & \textbf{Description} \\
\hline
\endhead
\hline \multicolumn{2}{r}{{Continued on next page}} \\
\endfoot
\hline
\endlastfoot

$\mathcal{G}$ & OptiGraph abstraction \\
$\mathcal{N}(\mathcal{G})$ & Set of nodes in graph $\mathcal{G}$ \\
$\mathcal{E}(\mathcal{G})$ & Set of edges in graph $\mathcal{G}$ \\
$x_n$ & Decision variables on node $n$ \\
$\mathcal{X}_n$ & Feasible set for decision variables on node $n$ \\
$f(\cdot)$ & Objective function \\
$g(\cdot)$ & Constraint function \\

$x_m$ & Master problem decision variables \\
$x_w$ & Subproblem variables for subproblem $w$ \\
$W$ & Set of subproblems \\
$c_m, c_w$ & Cost vectors for master and subproblems \\
$C_w, D_w, q_w$ & Constraint matrices and vectors \\
$\alpha_w$ & Cost-to-go variables \\
$\underline{\phi}_{m}^{i}$ & Master problem objective value at iteration $i$ \\
$\overline{\phi}_w^i$ & Subproblem objective value for subproblem $w$ at iteration $i$ \\
$\overline{x}_m^i$ & Master problem solution at iteration $i$ \\
$\lambda_w^i$ & Dual variables for subproblem $w$ at iteration $i$ \\
$UB^i$ & Upper bound at iteration $i$ \\
$LB^i$ & Lower bound at iteration $i$ \\

$\theta$ & Vector of uncertain/unknown parameters \\
$x$ & Vector of decision variables in mp problem \\
$J^*(\theta)$ & Optimal objective value as a function of $\theta$ \\
$\Theta$ & Compact polytope of bounded uncertain parameters \\
$z(\theta)$ & Parametric objective function in mp-LP \\
$c, H$ & Matrices and vectors in mp-LP objective function formulation \\
$A, b, F$ & Matrices and vectors inequality constrains in mp-LP formulation \\
$A_{eq}, b_{eq}, F_{eq}$ & Matrices and vectors for equality constraints in mp-LP \\
$A_{\theta}, b_{\theta}$ & Matrix and vector defining the parameter bounds space $\Theta$ \\
$v$ & Index for a critical region \\
$n_{CR}$ & Total number of critical regions \\
$x^*$ & Optimal solution of the mp problem \\
$A^v, b^v$ & Affine function matrices for solution in critical region $v$ \\
$E^v, f^v$ & Inequality matrices defining critical region $v$ \\

$\mathcal{P}$ & Set of processes \\
$\mathcal{J}$ & Set of chemicals \\
$\mathcal{T}$ & Set of planning periods \\
$\mathcal{K}$ & Set of uncertainty scenarios \\

$x_{pt}$ & Capacity expansion decisions for process $p$ in period $t$ \\
$y_{pt}$ & Binary decision to expand process $p$ in period $t$ \\
$q_{pt}$ & Cumulative capacity expansion for process $p$ in period $t$ \\

$p_{pt}$ & Production level of process $p$ in period $t$ \\
$b_{jt}$ & Amount of chemical $j$ purchased in period $t$ \\
$s_{jt}$ & Amount of chemical $j$ sold in period $t$ \\

$\alpha_{pt}, \beta_{pt}$ & Capital and variable cost coefficients for expansion \\
$\sigma_{pt}$ & Variable operating cost for process $p$ in period $t$ \\
$A_{jt}^k$ & Availability of chemical $j$ in period $t$ under scenario $k$ \\
$D_{jt}^k$ & Demand of chemical $j$ in period $t$ under scenario $k$ \\
$\gamma_{jt}^k$ & Purchase price of chemical $j$ in period $t$ under scenario $k$ \\
$\varphi_{jt}^k$ & Selling price of chemical $j$ in period $t$ under scenario $k$ \\
$\eta_{jp}$ & Stoichiometric consumption coefficient of chemical $j$ in process $p$ \\
$\mu_{jp}$ & Stoichiometric yield coefficient of chemical $j$ in process $p$ \\
$E_{pt}^L, E_{pt}^U$ & Lower and upper cumulative capacity bounds \\
$Q_{pt}^U$ & Per-period expansion limits \\
$\pi_k$ & Probability of scenario $k$ \\
$V(q, \xi)$ & Recourse value function for a given plan $q$ and realization $\xi$ \\
$\xi$ & Realization of uncertain parameters $(A_{jt}^k, D_{jt}^k, \gamma_{jt}^k, \varphi_{jt}^k)$ \\

\end{longtable}

\bibliography{BendersMP}

\end{document}